\documentclass[12pt]{amsart}
\usepackage{amsmath,amsthm,amssymb}
\usepackage{enumerate}
\usepackage{mathtools}
\usepackage{hyperref}
\mathtoolsset{showonlyrefs}

\newtheorem{theorem}{Theorem}
\newtheorem{lemma}[theorem]{Lemma}
\newtheorem{corollary}[theorem]{Corollary}
\newtheorem{proposition}[theorem]{Proposition}
\theoremstyle{definition}
\newtheorem{remark}[theorem]{Remark}
\newtheorem*{claim}{Claim}

\numberwithin{equation}{section}
\numberwithin{theorem}{section}
\allowdisplaybreaks

\newcommand{\ip}[2]{\langle #1,#2\rangle}
\DeclareMathOperator{\tr}{tr}
\DeclareMathOperator{\diam}{diam}
\newcommand{\e}{\varepsilon}

\author{Huy The Nguyen}
\address{School of Mathematical Sciences\\
	Queen Mary University of London\\
	Mile End Road\\
	London E1 4NS}
\email{h.nguyen@qmul.ac.uk}

\author{Artemis A. Vogiatzi}
\address{Copenhagen Centre for Geometry and Topology\\
	University of Copenhagen\\
	Universitetsparken 5, 2100 K{\o}benhavn, Danmark}
\email{aav@math.ku.dk}

\begin{document}
\title[]{Quadratic Convexity for Closed and Free Boundary High Codimension Mean Curvature Flow}

\begin{abstract}
We introduce quadratic convexity for mean curvature flow in high codimension. It combines quantitative convexity in the mean curvature direction with control of the full second fundamental form. Strict Andrews--Baker quadratic pinching implies quadratic convexity, while the new inequality also admits curvature tensors beyond the critical quadratic ratio.

For compact connected closed submanifolds of dimension $n\geq5$, we prove preservation of quadratic convexity and convergence in finite time to a round point. For compact connected submanifolds of dimension $n\geq4$ with spherical free boundary, we prove preservation and convergence to a round hemisphere in a uniquely determined $(n+1)$-plane. To our knowledge, this is the first such convergence theorem for spherical free boundary MCF in arbitrary codimension. The closed theorem requires additional quadratic pinching for $5\leq n\leq7$, and the free boundary theorem for $4\leq n\leq7$. For $n\geq8$, the necessary quadratic bound follows from quadratic convexity.

We also establish preservation of ordinary quadratic pinching $|A|^2\leq c|H|^2$ under spherical free boundary flow for $n\ge2$ and $\frac{9}{8n}\leq c\leq\frac{4}{3n}$, including both endpoints.
\end{abstract}
\maketitle
\setcounter{tocdepth}{1}
\tableofcontents

\section{Introduction}

Mean curvature flow (MCF) is the gradient descent flow of area functional. Given a smooth immersion $F_0\colon M^n\to\mathbb R^{n+m}$, a one parameter family of immersions $F(\cdot,t)$ satisfies

	\begin{align*}
\partial_t F= \Delta_{g(t)}F=H,\quad	F(\cdot,0)=F_0,
	\end{align*}
where $g(t)$ is the induced metric and $H=\tr_{g(t)}A$ is the mean curvature vector. For a closed immersion the first variation gives
	\begin{align*}
	\frac{d}{dt}\operatorname{Area}(M_t)&=-\int_{M_t}|H|^2\,d\mu_t.
	\end{align*}
which shows that the area decreases along the flow. The eventual shape of the submanifold near a singularity depends on what curvature information can be preserved during the flow.

In his seminal work, Huisken \cite{Huisken1984} showed that convexity is preserved by the flow. In fact he showed that a smooth, compact, uniformly convex hypersurface without boundary remains convex under MCF and contracts to a point in finite time. After rescaling, the hypersurfaces converge smoothly to a round sphere.

In higher codimension the second fundamental form takes values in a vector valued normal bundle. Its components in different normal directions interact in the curvature evolution, making preservation more difficult to establish. Andrews--Baker \cite[Theorem 1]{AndrewsBaker} discovered a quadratic condition on the norm of the second fundamental form that is preserved. In dimensions $n\ge4$, they showed that
	\begin{align*}
	|A|^2&\leq \frac{1}{n-1}|H|^2
	\end{align*}
is preserved by MCF of a compact connected submanifold without boundary whose mean curvature vector is initially nowhere zero. The flow contracts to a point in finite time, and the rescaled immersions converge to a round sphere in an $(n+1)$-plane.

More recently, Chow \cite{Chow2026} obtained convergence to a round point for closed flows in dimensions $n\ge10$ under tensorial curvature conditions extending classical scalar pinching.

The free boundary problem provides a natural setting in which to seek an analogue of these convergence results. Here the boundary of the evolving submanifold lies on a fixed supporting hypersurface, which it meets orthogonally. This condition eliminates the boundary term in the first variation of area, so the evolution remains its gradient flow. In codimension one, Stahl established existence and regularity for this problem \cite{Stahl1996} and proved that compact convex hypersurfaces with spherical free boundary contract to a point in finite time, with hemispherical singularity models \cite{StahlConvergence1996}. His result provides the free boundary counterpart of Huisken's theorem.

The subsequent development of the hypersurface theory has given a detailed understanding of curvature estimates and boundary singularities. Buckland \cite{Buckland} established a boundary monotonicity formula and studied singularity formation. Edelen \cite{EdelenConvexity2016} extended the Huisken--Sinestrari convexity estimates to free boundary flow and gave another proof of asymptotic umbilicity for convex flows with spherical boundary. More recently, Edelen--Haslhofer--Ivaki--Zhu \cite{Edelenetal} developed the regularity and structure theory for free boundary flow of mean convex hypersurfaces. In higher codimension, Evans--Lambert--Wood \cite{EvansLambertWood2022} proved convergence for equivariant Lagrangian discs under mixed Dirichlet and Neumann boundary conditions.

Our purpose is twofold: to introduce quadratic convexity as an extension of the Andrews--Baker condition for closed submanifolds, and to establish a corresponding convergence theorem with spherical free boundary in arbitrary codimension. The latter is a principal contribution of this paper. To our knowledge, this is the first theorem to establish convergence to a round hemisphere for submanifolds of arbitrary codimension satisfying a curvature pinching condition and meeting a fixed sphere orthogonally along their boundary. Under the initial curvature assumptions of Theorem~\ref{thm_fb-convergence}, in dimensions $n\geq4$, the flow contracts in finite time to a point on the supporting sphere. The rescaled immersions converge smoothly, including at the boundary, to a round hemisphere in a uniquely determined $(n+1)$-dimensional linear subspace.

Passing from hypersurfaces to higher codimension requires a curvature condition that controls the interaction between different normal components and is compatible with the free boundary geometry. Quadratic convexity combines quantitative convexity in the mean curvature direction with control of the remaining components of the second fundamental form. Its preservation uses both the interior evolution and the spherical conormal identities. Ordinary quadratic pinching illustrates why this boundary compatibility matters: although it is preserved in a substantial range, the boundary calculation imposes a lower bound on the pinching coefficient which excludes the critical Andrews--Baker coefficient in sufficiently high dimensions.

This additional requirement can be seen in the quadratic pinching calculation. For free boundary on a sphere of radius $R$, let $P=|A|^2-c|H|^2$ and denote the outward unit conormal by $\mu$. The boundary identity \eqref{eq_quadratic-boundary} gives
	\begin{align*}
	\frac R2\,\partial_\mu P&=-P-n\left|A_{\mu\mu}-\frac{3}{2n}H\right|^2 +\left(\frac9{4n}-2c\right)|H|^2.
	\end{align*}
Together with the interior evolution, this identity gives preservation of $P\leq0$ for $n\ge2$ and $\frac{9}{8n}\leq c\leq\frac{4}{3n}$, including both endpoints. The critical Andrews--Baker coefficient $\frac{1}{n-1}$ lies in this interval for $4\leq n\leq9$. In dimensions $n\ge10$ it falls below the lower endpoint, and the boundary identity no longer gives the sign required at $P=0$. The coarser bound with $c=\frac{4}{3n}$ remains available, but it permits the round cylinder $\mathbb S^{n-1}\times\mathbb R$.

For a hypersurface, weak convexity is equivalent to nonnegativity of the smallest principal curvature, with the normal chosen in the mean curvature direction wherever $|H|>0$. In high codimension, the mean curvature vector also selects a normal direction. Set $\nu=H/|H|$ and decompose $A=h\nu+A^-$, where $h=\langle A,\nu\rangle$. If $\lambda_1$ is the smallest eigenvalue of $h$, then $\lambda_1\ge0$ expresses convexity in the mean curvature direction.

Convexity in the principal normal direction alone does not force convergence to a round point. For $n\ge4$ and $0<r<1$, the compact products $M_r=\mathbb S^{n-1}(r)\times\mathbb S^1(1)\subset\mathbb R^{n+2}$ satisfy
	\begin{align*}
	|H|^2&=\frac{(n-1)^2}{r^2}+1, \qquad |A|^2=\frac{n-1}{r^2}+1,\\
	\lambda_1&=\frac1{|H|}, \qquad \lambda_2=\cdots=\lambda_n=\frac{n-1}{r^2|H|}.
	\end{align*}
Thus $h$ is positive definite, but their MCF is
	\begin{align*}
	M_r(t)&=\mathbb S^{n-1}\!\left(\sqrt{r^2-2(n-1)t}\right) \times\mathbb S^1\!\left(\sqrt{1-2t}\right), \qquad 0\leq t<\frac{r^2}{2(n-1)}.
	\end{align*}
At the first singular time it collapses onto a circle of positive radius, rather than to a round point.

The spherical boundary identities give a reason to use $\lambda_1$ directly. At a boundary point where $|H|>0$, $h\ge0$ and $\lambda_1=0$, the vanishing of $h_{a\mu}$ allows us to choose a unit null eigenvector $e$ tangent to $\partial M$ or equal to $\mu$. Since $\nabla_\mu^\perp\nu=0$, the identities of Section \ref{sec:quadratic-pinching} give
	\begin{align*}
	R(\nabla_\mu h)(e,e) = \begin{cases} h_{\mu\mu}\ge0,&e\in T\partial M,\\
	2|H|>0,&e=\mu. \end{cases}
	\end{align*}
These are the conormal signs compatible with preservation of convexity. The interior evolution of $h$ also couples to the other normal components of $A$, so we combine the least eigenvalue with additional curvature control. For $\kappa>0$ and $\e\ge0$, define the quadratic convexity quantity
	\begin{align}\label{eqn_F}
	F_\e&=|H|\lambda_1-\frac{n-1}{2}\bigl(\lambda_1^2+|A^-|^2\bigr)+\kappa\left[\left(\frac1{n-1}-\e\right)|H|^2-|A|^2\right].
	\end{align}
We refer to $F_\e>0$ as quadratic convexity. We interpret this inequality as positivity of $|H|^2F_\e$, extended continuously by zero where $H=0$.

The algebraic estimates in Section \ref{sec:notation} relate quadratic convexity to critical quadratic pinching and convexity in the mean curvature direction. At points where $|H|>0$,
	\begin{align*}
	|A|^2&\leq\frac{|H|^2}{n-1} \quad\Longrightarrow\quad F_0\ge0 \quad\Longrightarrow\quad h\ge0.
	\end{align*}
The condition $F_0>0$ can hold beyond the critical quadratic ratio. The new preservation result therefore enlarges the class of curvature tensors permitted by the Andrews--Baker inequality while retaining quantitative convexity in the mean curvature direction. It still imposes more than convexity, even when $A^-=0$. For example, positive definite hypersurface tensors of trace one approaching $\operatorname{diag}(1,0,\ldots,0)$ have $F_0<0$ for $n\ge4$ once they are sufficiently close to that tensor.

The same products give a precise sharpness statement for closed flows. Since $\lambda_2=\cdots=\lambda_n$, the sum of squares involving the remaining eigenvalues in the trace decomposition in Section \ref{sec:notation} vanishes. Consequently, on $M_r$,
	\begin{align*}
	F_0&=-(n-2)\left(\frac{\kappa}{n-1}+\frac12\right) +\frac{r^2}{2\bigl((n-1)^2+r^2\bigr)}<0,\\*
	\frac{F_0}{|H|^2}&\longrightarrow0\quad\text{ from below as }r\to0.
	\end{align*}
Thus, for every $\delta>0$, there are compact initial submanifolds that are strictly convex in the principal normal direction and satisfy $F_0>-\delta|H|^2$ whose flows do not converge to round points. The sign condition is therefore sharp with respect to these perturbations. Moreover,
	\begin{align*}
	\frac{|A|^2}{|H|^2}&=\frac{n-1+r^2}{(n-1)^2+r^2} \longrightarrow\frac1{n-1}.
	\end{align*}
For $5\leq n\le7$, these examples also satisfy the additional quadratic hypothesis of Theorem \ref{thm_closed-convergence} when $r$ is sufficiently small; in dimension four they lie outside that hypothesis. Since $F_\e=F_0-\kappa\e|H|^2$, the inequality $F_\e>0$ gives a quantitative lower bound for $F_0/|H|^2$ when $\e>0$ is fixed. This example establishes sharpness of the sign condition for closed flows, without asserting optimality of $\kappa$ or a fixed $\e$, or sharpness for flows with spherical free boundary.

On a compact initial submanifold, strict Andrews--Baker pinching implies $F_\e>0$ for sufficiently small $\e>0$ (Lemma \ref{lem_initial-positivity}). For fixed $\e>0$, quadratic convexity gives the quantitative bound
	\begin{align*}
	F_\e&\ge0 \quad\Longrightarrow\quad h\geq \frac{\kappa\e}{1+2\kappa/(n-1)}\,|H|g \qquad (|H|>0).
	\end{align*}
This bound excludes the round cylinder, which has $\lambda_1=0$. In dimensions $n\ge8$, the same condition also gives $|A|^2<\frac{4}{3n}|H|^2$ for $\kappa=14/5$, providing the quadratic control used in the convergence proofs. For $4\leq n\leq7$, we assume this quadratic bound separately in the theorem for flows with free boundary.

In both theorems we work with codimension $m\ge2$ and $\kappa=14/5$ in the definition of $F_\e$.

\begin{theorem}[Convergence to a round point for closed flows]
\label{thm_closed-convergence}
Let $F:M^n\times[0,T)\to\mathbb R^{n+m}$ be the maximal smooth MCF of a compact connected manifold without boundary, with $n\ge5$.
\begin{enumerate}
\item For $n \geq 8$ assume $F_\e>0$ initially for a fixed $\e>0$, sufficiently small.

\item For $5\leq n\leq7$ assume $F_\e>0$ initially for a fixed $\e>0$, sufficiently small and in addition that $|A|^2<c_n|H|^2$ initially, where
	\begin{align*}
	c_n&=\frac{3(n+1)}{2n(n+2)}, \quad\text{ for } 5\leq n\leq7.
	\end{align*}
\end{enumerate}
Then $F_\e>0$ is preserved, $T<\infty$, and the flow contracts to a point $q\in\mathbb R^{n+m}$. After reparametrisation,
	\begin{align*}
	\frac{F(\cdot,t)-q}{\sqrt{2n(T-t)}}
	\end{align*}
converges in $C^\infty$ to a parametrisation of the unit round $n$-sphere in a fixed $(n+1)$-dimensional linear subspace.
\end{theorem}
\begin{remark}
\begin{enumerate}
\item In dimensions $2\le n\le4$, the additional quadratic pinching coincides with the strict Andrews--Baker condition. Hence, convergence to a round point already follows from \cite[Theorem 1]{AndrewsBaker} and we do not include this dimension as a separate improvement.

\item The theorem also applies in codimension one. However, in this case $F_\e>0$ implies strict convexity, so convergence to a round point already follows from Huisken's theorem \cite{Huisken1984}. The novelty of our result therefore lies in higher codimension.
\end{enumerate}
\end{remark}

\begin{theorem}[Convergence to a round hemisphere]
\label{thm_fb-convergence}
Let $F:M^n\times[0,T)\to\overline B_R^{n+m}$, $m\ge2$, be the maximal smooth MCF of a compact connected manifold with nonempty orthogonal free boundary on $\mathbb S_R^{n+m-1}$, with $n\ge4$. Assume $F_\e>0$ initially for some sufficiently small $\e>0$. Assume in addition that $|A|^2<\frac{4}{3n}|H|^2$ holds initially, for $4\leq n\le7$. Then $F_\e>0$ is preserved, $T<\infty$, and the flow contracts to a point $q\in\mathbb S_R^{n+m-1}$. After reparametrisation,
	\begin{align*}
	\frac{F(\cdot,t)-q}{\sqrt{2n(T-t)}}
	\end{align*}
converges in $C^\infty$ up to the free boundary to a parametrisation of the unit round hemisphere
	\begin{align*}
	\{x\in V: |x|&=1,\ \langle x,q\rangle\leq0\},
	\end{align*}
where $V\subset\mathbb R^{n+m}$ is a fixed $(n+1)$-dimensional linear subspace containing $q/R$.
\end{theorem}
\begin{remark}
The one-dimensional analogue of Theorem \ref{thm_fb-convergence} was established in \cite{FBHCCSF} for high-codimension curve shortening flow with free boundary. Under a reflected entropy bound below $2$, finite-time singularities have shrinking round semicircles as blow-up models.
\end{remark}

In Section \ref{sec:notation} we fix the notation and establish the algebraic properties of $F_\e$. We show how strict Andrews--Baker pinching implies quadratic convexity and how quadratic convexity controls the least principal curvature in the mean curvature direction. In dimensions $n\ge8$ it also supplies the quadratic bound needed for the curvature estimates. In Section \ref{sec:interior} we estimate the reaction and gradient terms at a possible first zero of $F_\e$. The resulting interior differential inequality is used in both preservation arguments. For closed flows, we combine it with existing curvature estimates to prove convergence to a round point.

In Section \ref{sec:quadratic-pinching} we derive the conormal curvature identities for a spherical free boundary. These give preservation of ordinary quadratic pinching in the stated range and an additive gap under strict initial pinching, which keeps $|H|$ bounded away from zero. In Section \ref{sec:planarity} we use the sign of the conormal derivative of the planarity quotient to extend the interior planarity estimate up to the free boundary. This ensures that the component of $A$ orthogonal to $H$ disappears in curvature blowups. In Section \ref{sec:gradient-estimates} we construct corrections to the curvature derivatives that allow the maximum principle to give bounds for $\nabla A$ and $\nabla^2 A$ up to the boundary.

In Section \ref{sec:boundary-proof} we compute the conormal derivative of an upper support for $F_\e$ at a possible first zero. Its sign excludes boundary contact and completes the preservation proof. In Section \ref{sec:fb-convergence} we use the derivative estimates to obtain smooth blowup limits up to the free boundary and prove that every limit at the extinction scale is a round hemisphere. Finally, the quantitative planarity estimate controls the change of the ambient $(n+1)$-plane near extinction. This makes the limiting plane independent of the blowup sequence and gives smooth convergence of the whole rescaled flow.

Appendix \ref{app:fb-compactness} establishes the compactness and continuation result used in the convergence proof, including the adaptations of graphical compactness and parabolic regularity to the free boundary.

\section*{Declarations}
\subsection*{Funding} The first author was financially supported by the EPSRC through the grant ``Geometric Flows and the Dynamics of Phase Transitions" (EP/Y017862/1). The second author is financially supported by the Danish National Research Foundation through the Copenhagen Centre for Geometry and Topology (DNRF151).

\subsection*{AI disclosure}
The authors used LLMs to assist with checking mathematical arguments and calculations, identifying relevant literature and revising the exposition. The authors take full responsibility for the mathematical content, references and final manuscript.

\section{Notation and Algebraic Pinching Estimates}\label{sec:notation}

Let $F:M^n\times[0,T)\to\mathbb R^{n+m}$ be a smooth MCF, and let $F_t=F(\cdot,t)$ and $M_t=F_t(M)$. We write $g$ for the induced metric, $D$ for the Euclidean connection, and $\nabla$ and $\nabla^\perp$ for the induced tangent and normal connections. Thus $D_XY=\nabla_XY+A(X,Y)$, $H=\tr_g A$, and $\partial_tg_{ij}=-2\langle H,A_{ij}\rangle$. We use $\nabla$ also for the connection on covariant tensors with values in the normal bundle. In orthonormal frames, tangent indices range from $1$ to $n$ and normal indices from $1$ to $m$. The symmetric matrix $A_\alpha$ has entries $A_{ij\alpha}=\langle A(e_i,e_j),\nu_\alpha\rangle$. Matrix products and inner products contract the tangent indices.

In the operator $\mathcal L=\nabla_{\partial_t}-\Delta$, the time connection is induced by ambient orthogonal projection and preserves the metric. Here the Laplacian is the rough Laplacian. On functions, $\mathcal L=\partial_t-\Delta$. We use the evolution conventions of \cite[Section 2]{LN}. By Euclidean Codazzi, $(\nabla_k A)_{ij}$ is symmetric in the three tangent indices.

In the decomposition along the principal normal, $\tr_g h=|H|$, $\tr_g A^-=0$, and $|A|^2=|h|^2+|A^-|^2$. Set $\mathring h=h-|H|g/n$ and order the eigenvalues $\lambda_1\le\cdots\le\lambda_n$ of $h$. In a normal frame with $\nu_1=\nu$, we have $A_1=h$ and $A_\alpha=A_\alpha^-$ for $\alpha>1$. For a local smooth unit tangent field $e$, write $F_{\e,e}$ for \eqref{eqn_F} with $\lambda_1$ replaced by $h(e,e)$.

For spherical free boundary flow, $F(\operatorname{int}M,t)\subset B_R^{n+m}$ and $F(\partial M,t)\subset\mathbb S_R^{n+m-1}$. The orthogonality condition is $\mu=F/R$ along $\partial M_t$, where $\mu$ is the outward unit conormal, see \cite{Buckland,Stahl1996}. In a boundary frame $e_n=\mu$ and the indices $a,b,c$ range from $1$ to $n-1$. Since $D_{e_a}\mu=e_a/R$, we have $\nabla_{e_a}\mu=e_a/R$ and $A_{a\mu}=0$. Under the parabolic rescaling $\widetilde F(p,\tau)=R^{-1}F(p,R^2\tau)$, the norms of the covariant derivatives of $A$ scale by $|\widetilde{\nabla^k A}|=R^{k+1}|\nabla^k A|$.

\subsection{Algebraic consequences of the pinching}

Since $\sum_{i=2}^n\lambda_i=|H|-\lambda_1$ and $|A|^2=|A^-|^2+\sum_{i=1}^n\lambda_i^2$, expanding the squares gives
	\begin{align*}
	\sum_{i=2}^n\left(\lambda_i-\frac{|H|-\lambda_1}{n-1}\right)^2&=|A|^2-|A^-|^2-\lambda_1^2 -\frac{(|H|-\lambda_1)^2}{n-1}.
	\end{align*}
Solving for $|A^-|^2$ and substituting into the term $-(n-1)|A^-|^2/2$ in \eqref{eqn_F}, we use
	\begin{align*}
	|H|\lambda_1+\frac12(|H|-\lambda_1)^2&=\frac12|H|^2+\frac12\lambda_1^2
	\end{align*}
to obtain
	\begin{align*}
	F_\e={}&\left(\kappa+\frac{n-1}{2}\right) \left(\frac{|H|^2}{n-1}-|A|^2\right) +\frac12\lambda_1^2-\kappa\e|H|^2\\
	&+\frac{n-1}{2}\sum_{i=2}^n \left(\lambda_i-\frac{|H|-\lambda_1}{n-1}\right)^2.
	\end{align*}
Setting $\e=0$ gives the decomposition of $F_0$. The following lemma makes precise how strict Andrews--Baker pinching on a compact initial submanifold implies $F_\e>0$ for sufficiently small $\e>0$.

\begin{lemma}[Quadratic pinching implies quadratic convexity]
\label{lem_initial-positivity}
Let $M_0^n$, $n\ge2$, be a compact immersed submanifold, and fix $\kappa>0$. Suppose
	\begin{align*}
	|A|^2&<\frac{|H|^2}{n-1}, \qquad\text{ on }M_0.
	\end{align*}
Then, there exists $\e_0>0$, such that $F_\e>0$, on $M_0$, for every $0<\e\le\e_0$. More generally, if only
	\begin{align*}
	|A|^2&\le\frac{|H|^2}{n-1}
	\end{align*}
is assumed, such an $\e_0$ exists if and only if the principal second fundamental form $h$ is defined and positive definite everywhere on $M_0$.
\end{lemma}

\begin{proof}
Under strict quadratic pinching, compactness gives
	\begin{align*}
	b&=\min_{M_0}\left(\frac{|H|^2}{n-1}-|A|^2\right)>0.
	\end{align*}
Since $|A|^2\ge0$, we have $|H|^2\ge(n-1)b$, so the principal normal is defined everywhere. The preceding decomposition gives $F_0\ge(\kappa+(n-1)/2)b>0$. By compactness,
	\begin{align*}
	\min_{M_0}\frac{F_0}{|H|^2}&>0.
	\end{align*}
The identity
	\begin{align*}
	\frac{F_\e}{|H|^2}&=\frac{F_0}{|H|^2}-\kappa\e
	\end{align*}
therefore proves the assertion for sufficiently small $\e_0>0$.

Now assume the nonstrict quadratic pinching. If $H=0$ at any point, the strict inequality $F_\e>0$ fails there by definition. We may therefore restrict to the case where $H$ is nowhere zero. Since $\operatorname{tr}h=|H|$, Cauchy-Schwarz gives
	\begin{align*}
	0&\le\frac{|H|^2}{n-1}-|A|^2\le\frac{|H|^2}{n-1} -\lambda_1^2 -\frac{(|H|-\lambda_1)^2}{n-1} -|A^-|^2=\frac{\lambda_1(2|H|-n\lambda_1)}{n-1} -|A^-|^2.
	\end{align*}
Because $\lambda_1\leq |H|/n$, the factor $2|H|-n\lambda_1$ is positive. Hence $\lambda_1\ge0$.

If $\lambda_1=0$ at a point, equality is forced throughout the preceding estimate. Consequently,
	\begin{align*}
	A^-&=0, \qquad \lambda_2=\cdots=\lambda_n=\frac{|H|}{n-1}, \qquad |A|^2=\frac{|H|^2}{n-1}.
	\end{align*}
At that point,
	\begin{align*}
	F_0&=0, \qquad F_\e=-\kappa\e|H|^2<0
	\end{align*}
for every $\e>0$. Thus positive definiteness of $h$ is necessary.

Conversely, if $h$ is positive definite everywhere, the preceding decomposition and the nonstrict quadratic pinching give
	\begin{align*}
	F_0&\ge\frac12\lambda_1^2>0.
	\end{align*}
The same compactness argument then provides a uniform $\e_0>0$, for which $F_\e>0$.
\end{proof}

\begin{lemma}[Principal convexity from quadratic convexity]
\label{lem_modified-principal-convexity}
At a point where $|H|>0$, suppose $n\ge2$, $\kappa>0$, $\e>0$ and $F_\e\ge0$. Then
	\begin{align}\label{rangeoflambda1}
	0&<\frac{\kappa\e}{1+2\kappa/(n-1)}|H| \le\lambda_1\le\frac{|H|}{n}.
	\end{align}
If $e$ is a local smooth unit vector field agreeing at the point with a minimum eigenvector of $h$, then $F_{\e,e}$ is a smooth upper support for $F_\e$ at that point.
\end{lemma}
\begin{proof}
Since the other $n-1$ eigenvalues sum to $|H|-\lambda_1$, Cauchy-Schwarz gives
	\begin{align*}
	|A|^2&=|h|^2+|A^-|^2 \ge\lambda_1^2+\frac{(|H|-\lambda_1)^2}{n-1}+|A^-|^2.
	\end{align*}
Substituting this into \eqref{eqn_F}, we obtain
	\begin{align*}
	F_\e\le{}&\left(1+\frac{2\kappa}{n-1}\right)|H|\lambda_1 -\left(\frac{n-1}{2}+\frac{n\kappa}{n-1}\right)\lambda_1^2-\left(\kappa+\frac{n-1}{2}\right)|A^-|^2 -\kappa\e|H|^2.
	\end{align*}
Using $F_\e\ge0$ and discarding the nonpositive square terms proves the lower bound in \eqref{rangeoflambda1}. The upper bound follows from $\lambda_1\le(\operatorname{tr}h)/n=|H|/n$.

In particular,
	\begin{align*}
	|H|-(n-1)\lambda_1&\ge\frac{|H|}{n}>0.
	\end{align*}
Thus $s\mapsto |H|s-(n-1)s^2/2$ is increasing near $s=\lambda_1$. Since $h(e,e)\ge\lambda_1$ nearby and equality holds at the point, continuity gives $F_{\e,e}\geq F_\e$ in a neighbourhood, with equality at the point. This proves the upper support assertion.
\end{proof}

\begin{lemma}[Quadratic pinching from quadratic convexity]
\label{lem_modified-quadratic-pinching}
At a point where $|H|>0$, let $\kappa>0$ and $\e>0$. If the quantity in \eqref{eqn_F} satisfies $F_\e\ge0$, then
	\begin{align}\label{eq_modified-quadratic-pinching}
	|A|^2&\le \left(\frac1{n-1}+\frac{n+1}{2\kappa n^2}-\e\right)|H|^2.
	\end{align}
In particular, if $n\ge8$ and $\kappa=14/5$, then
	\begin{align}\label{eq_modified-strict-coarse}
	|A|^2&<\frac4{3n}|H|^2.
	\end{align}
\end{lemma}
\begin{proof}
Since $\lambda_1$ is the least eigenvalue and $\operatorname{tr}h=|H|$, we have $\lambda_1\le|H|/n$. The function $s\mapsto |H|s-(n-1)s^2/2$ is increasing for $s\le|H|/n$. Consequently
	\begin{align*}
	|H|\lambda_1-\frac{n-1}{2}\lambda_1^2&\leq \frac{n+1}{2n^2}|H|^2.
	\end{align*}
Rearranging $F_\e\ge0$ and using this inequality gives
	\begin{align*}
	\kappa|A|^2&\leq \kappa\left(\frac1{n-1}-\e\right)|H|^2 +|H|\lambda_1-\frac{n-1}{2}\lambda_1^2 -\frac{n-1}{2}|A^-|^2\\
	&\leq \left[\kappa\left(\frac1{n-1}-\e\right) +\frac{n+1}{2n^2}\right]|H|^2.
	\end{align*}
Dividing by $\kappa$ proves \eqref{eq_modified-quadratic-pinching}. For $\kappa=14/5$,
	\begin{align*}
	\frac4{3n}-\frac1{n-1}-\frac{5(n+1)}{28n^2}&=\frac{13n^2-112n+15}{84n^2(n-1)}.
	\end{align*}
The numerator is positive and increasing for $n\ge9$, which proves \eqref{eq_modified-strict-coarse} in those dimensions.

For $n=8$, normalise $|H|=1$ and set $s=\lambda_1$. The positive semidefinite tensor $h-sg$ has trace $1-8s$, so $|h|^2\le\Phi(s):=1/8+56(1/8-s)^2$. If $|A|^2\ge1/6$, then $|A^-|^2\ge(1/6-\Phi(s))_+$, where $r_+=\max\{r,0\}$. Hence
	\begin{align*}
	F_\e&\leq s-\frac72s^2-\frac72(1/6-\Phi(s))_+-\frac1{15}-\kappa\e.
	\end{align*}
The first three terms are increasing for $s\leq s_*=(1-1/\sqrt{21})/8$, and form a convex quadratic on $[s_*,1/8]$. Their values at $s_*$ and $1/8$ are, respectively, $(13/3-1/\sqrt{21})/64$ and $-29/384$, both less than $1/15$. Thus $F_\e<0$, a contradiction.
\end{proof}

\section{Curvature Preservation}\label{sec:interior}

\subsection{Interior preservation}
The following interior inequality is the main curvature calculation in preserving quadratic convexity. At a first zero, the combined reaction and gradient terms control both $A^-$ and the variation of the principal normal $\nu$. We denote these contributions by $\mathcal R(F_\e)$ and $\mathcal G(F_\e)$. The estimate applies independently of the boundary condition.
\begin{lemma}[Interior differential inequality]\label{lem_interior}
Let $n\ge4$ and fix $\kappa=14/5$. At an interior zero of \eqref{eqn_F}, suppose $|H|>0$ and $|A|^2\leq \frac{4}{3n} |H|^2$. Choose a unit minimum eigenvector $e$ of $h$ and extend it with $\nabla e=0$ and $\nabla_{\partial_t}e=0$ at the point, using the time connection that preserves the metric. For sufficiently small $\e>0$ we have
	\begin{align}\label{finalestimate}
	(\partial_t-\Delta)F_{\e,e}&\geq c_1|H|^2|A^-|^2 +c_2|H|^2|\nabla^\perp\nu|^2\ge0,
	\end{align}
for constants $c_1,c_2>0$ depending only on $n$. The smallness restriction on $\e$ is needed only for $n=4$. For $n\ge8$, the quadratic pinching at the zero follows from \eqref{eq_modified-strict-coarse}.
\end{lemma}

\subsection{Estimates at a zero of $F_\e$ and choice of $\kappa$}
Let $\kappa>0$. At $F_\e=0$, \eqref{eqn_F} becomes
	\begin{align}\label{eqn_contact}
	\kappa\left(\frac{|H|^2}{n-1}-|A|^2\right) -\frac{n-1}{2}|A^-|^2&=\kappa\e|H|^2-|H|\lambda_1+\frac{n-1}{2}\lambda_1^2.
	\end{align}
Recall \eqref{rangeoflambda1}. Since $s-(n-1)s^2/2$ is increasing on $[0,1/n]$, substitution into \eqref{eqn_contact} gives
	\begin{align}\label{eqn_contact-smallness}
	\frac{|A^-|^2}{|H|^2}&\leq U_n-\beta_n\frac{|\mathring h|^2}{|H|^2}, \quad U_n=\frac1{n(n-1)} +\frac1{2n^2((n-1)/2+\kappa)},\\
	\beta_n&=\frac{\kappa}{(n-1)/2+\kappa}.
	\end{align}
Set $x=\frac{|\mathring h|^2}{|H|^2}$. Multiplying \eqref{eqn_contact-smallness} by $x$ gives
	\begin{align*}
	\frac{|\mathring h|^2|A^-|^2}{|H|^4}&\leq U_nx-\beta_nx^2.
	\end{align*}
The quadratic on the right hand side attains the maximum at $x=U_n/(2\beta_n)$. Hence
	\begin{align}\label{eqn_contact-product}
	\frac{|\mathring h|^2|A^-|^2}{|H|^4}&\leq U_n\frac{U_n}{2\beta_n}-\beta_n\left(\frac{U_n}{2\beta_n}\right)^2 =\frac{U_n^2}{4\beta_n}.
	\end{align}
For the gradient estimate below it is convenient to require
	\begin{align*}
	A_n(\kappa):&=\frac{2\kappa}{n+2} -\frac{\kappa+1/2}{n-1}\ge\frac1n.
	\end{align*}
Let
	\begin{align}\label{eqn_kappa-choice}
	\kappa&=\frac{14}{5}.
	\end{align}
For $n\ge10$, this choice satisfies the coefficient inequality required below. The remaining dimensions $4\leq n\le9$ are treated separately in the completion of the proof. We next derive the required bounds directly from $F_\e\ge0$. Since
	\begin{align*}
	|A|^2&=\frac{|H|^2}{n}+|\mathring h|^2+|A^-|^2,
	\end{align*}
rearranging \eqref{eqn_F} gives
	\begin{align*}
	\kappa|\mathring h|^2 +\left(\kappa+\frac{n-1}{2}\right)|A^-|^2&\leq |H|\lambda_1-\frac{n-1}{2}\lambda_1^2 +\kappa\left(\frac1{n(n-1)}-\e\right)|H|^2\\
	&\leq \left[\frac{n+1}{2n^2} +\frac{\kappa}{n(n-1)}-\kappa\e\right]|H|^2.
	\end{align*}
Here the last inequality follows from $\lambda_1\leq |H|/n$. Therefore
	\begin{align}\label{eqn_F-tracefree}
	|\mathring h|^2 +\left(1+\frac{n-1}{2\kappa}\right)|A^-|^2&\leq \left(\frac1{n(n-1)} +\frac{n+1}{2\kappa n^2}-\e\right)|H|^2.
	\end{align}
For $\kappa=14/5$ and $n\ge10$,
	\begin{align*}
	\frac1{3n} -\frac1{n(n-1)} -\frac{5(n+1)}{28n^2}&=\frac{13n^2-112n+15}{84n^2(n-1)}>0,
	\end{align*}
because the numerator equals $195+148(n-10)+13(n-10)^2$. In particular, \eqref{eqn_F-tracefree} implies
	\begin{align*}
	|\mathring h|^2&\le\frac{|H|^2}{3n}.
	\end{align*}
For this choice of $\kappa$, the constants in \eqref{eqn_contact-smallness} become
	\begin{align*}
	U_n&=\frac1{n(n-1)}+\frac5{n^2(5n+23)},\qquad \beta_n=\frac{28}{5n+23}.
	\end{align*}
For $n\ge10$,
	\begin{align}\label{eqn_contact-uniform}
	U_n&\le\frac{17}{16n(n-1)},\qquad \frac1{4\beta_n}=\frac{5n+23}{112}.
	\end{align}
The first inequality is equivalent to $5n^2-57n+80\ge0$.

\subsection{Reaction terms}
\begin{lemma}[Estimate for the reaction terms]
Let $n\ge10$, $\kappa=14/5$ and $\e>0$. At an interior zero of $F_\e$ with $|H|>0$, we have
	\begin{align}\label{finalestimateR}
	\mathcal R(F_\e)&\geq c_1|H|^2|A^-|^2,
	\end{align}
for some $c_1>0$.
\end{lemma}
\begin{proof}
The evolution equations for $A$ and $H$ in \cite[Section 2]{LN} give the following reaction terms, where $\nu_1=\nu$, $A_1=h$, and $A_\alpha=A^-_\alpha$ for $\alpha>1$:
	\begin{align*}
	\mathcal R(A)_\alpha&=\sum_\beta\left[ \ip{A_\alpha}{A_\beta}A_\beta+A_\beta^2A_\alpha +A_\alpha A_\beta^2-2A_\beta A_\alpha A_\beta\right].
	\end{align*}
In particular,
	\begin{align}\label{normevolutions}
	\mathcal R(|H|)&=|h|^2|H|,\\
	\mathcal{R}(\nu)&=\sum_{\alpha>1}\left\langle h, A_\alpha^-\right\rangle \nu_\alpha\nonumber\\
	\mathcal{R}(|A^-|^2)&=2 \sum_{\alpha>1}\left|\left[h, A_\alpha^-\right]\right|^2+2 \sum_{\alpha, \beta>1}\left(\left\langle A_\alpha^-, A_\beta^-\right\rangle^2+\left|\left[A_\alpha^-, A_\beta^-\right]\right|^2\right).\nonumber
	\end{align}
Also setting $\alpha =1$ we have
	\begin{align*}
	\mathcal R(h)&=|h|^2h+\mathcal T,
	\end{align*}
where
	\begin{align*}
	\mathcal T&=2\sum_{\alpha>1}\ip{h}{A^-_\alpha}A^-_\alpha +\sum_{\alpha>1}\bigl((A^-_\alpha)^2h+h(A^-_\alpha)^2-2A^-_\alpha hA^-_\alpha\bigr).
	\end{align*}
Here differentiation of $H/|H|$ is included in $\mathcal R(h)$. Diagonalise $h$ and let $e=e_1$ be a minimum eigenvector. Then
	\begin{align*}
	\mathcal T_{11}&=2\sum_{\alpha>1} \ip{h}{A^-_\alpha}(A^-_\alpha)_{11}+2\sum_{\alpha>1,j}(\lambda_1-\lambda_j)(A^-_\alpha)_{1j}^2.
	\end{align*}
\begin{claim}
	\begin{align}\label{eqn_T-lower}
	\mathcal T_{11}&\geq |A^-|^2\left(\lambda_1-\frac{|H|}{n}-\sqrt{\frac{n-1}{n}}|\mathring h|\right).
	\end{align}
\end{claim}
\begin{proof}
For the diagonal part, let $x_i=(A^-_\alpha)_{ii}$ and write
	\begin{align*}
	\mu_i&=\lambda_i-\frac{|H|}{n}, \qquad \sum_i x_i=\sum_i\mu_i=0.
	\end{align*}
Then
	\begin{align*}
	2\langle h,A^-_\alpha\rangle(A^-_\alpha)_{11}&=2x_1\sum_i\mu_i x_i.
	\end{align*}
Writing $u=e_1-\frac1n(1,\ldots,1)$, we have $x_1=\langle u,x\rangle$, $\mu_1=\langle u,\mathring{h}\rangle$, $|u|=\sqrt{\frac{n-1}{n}}$. Hence
	\begin{align*}
	2x_1\sum_i\mu_i x_i&=2\langle u,x\rangle\langle\mathring{h},x\rangle\ge\bigl(\langle u,\mathring{h}\rangle-|u||\mathring{h}|\bigr)|x|^2=\left(\lambda_1-\frac{|H|}{n}-\sqrt{\frac{n-1}{n}}|\mathring h|\right)|x|^2.
	\end{align*}
For the terms from entries off the diagonal,
	\begin{align*}
	\lambda_1-\lambda_j&\ge\lambda_1-\lambda_n \ge\lambda_1-\frac{|H|}{n} -\sqrt{\frac{n-1}{n}}|\mathring h|,
	\end{align*}
since
	\begin{align*}
	\lambda_n-\frac{|H|}{n}&\le\sqrt{\frac{n-1}{n}}|\mathring h|.
	\end{align*}
Therefore, summing over $\alpha>1$,
	\begin{align}
	\mathcal T_{11}&\geq |A^-|^2\left(\lambda_1-\frac{|H|}{n}-\sqrt{\frac{n-1}{n}}|\mathring h|\right).
	\end{align}
which is \eqref{eqn_T-lower}.
\end{proof}
Furthermore, diagonalisation and Cauchy-Schwarz give
	\begin{align*}
	\sum_{\alpha>1}\left(\ip{h}{A^-_\alpha}^2+|[ h,A^-_\alpha]|^2\right)&\le2|\mathring h|^2|A^-|^2,\\
	\sum_{\alpha,\beta>1}\ip{A^-_\alpha}{A^-_\beta}^2&\le|A^-|^4,\\
	\sum_{\alpha,\beta>1}|[ A^-_\alpha,A^-_\beta]|^2&\le2|A^-|^4.
	\end{align*}
Substituting in \eqref{normevolutions} yields
	\begin{align}\label{eqn_norm-reactions}
	\mathcal R(|A^-|^2)&\le6|A^-|^4+4|\mathring h|^2|A^-|^2,
	\end{align}

\begin{claim}
	\begin{align}\label{eqn_norm-reactions2}
	\mathcal R\left(\frac{|H|^2}{n-1}-|A|^2\right)&\geq 2\bigl(|h|^2+3|A^-|^2\bigr)\left(\frac{|H|^2}{n-1}-|A|^2\right)\\
	&+2\left(\frac4n-\frac3{n-1}\right)|H|^2|A^-|^2.
	\end{align}
\end{claim}
Inequality \eqref{eqn_norm-reactions2} holds regardless of the sign of $|H|^2/(n-1)-|A|^2$.
\begin{proof}
To see this explicitly, we calculate
	\begin{align*}
	\mathcal R\left(\frac{|H|^2}{n-1}-|A|^2\right)&-2\bigl(|h|^2+3|A^-|^2\bigr)\left(\frac{|H|^2}{n-1}-|A|^2\right)-2\left(\frac4n-\frac3{n-1}\right)|H|^2|A^-|^2\\
	&=4\Bigg[2|\mathring h|^2|A^-|^2-\sum_{\alpha>1}\left(\ip{h}{A^-_\alpha}^2+|[h,A^-_\alpha]|^2\right)\Bigg]\\
	&+2\Bigg[3|A^-|^4-\sum_{\alpha,\beta>1} \left(\ip{A^-_\alpha}{A^-_\beta}^2 +|[A^-_\alpha,A^-_\beta]|^2\right)\Bigg]\ge0.
	\end{align*}
Each bracket is nonnegative by the matrix inequalities above, with no assumption on $|H|^2/(n-1)-|A|^2$.
\end{proof}
In the reaction of $F_\e$, the term proportional to $(|H|^2/(n-1)-|A|^2)|A^-|^2$ is retained. At $F_\e=0$ it combines with the quartic term by the exact identity
	\begin{align*}
	6|A^-|^2\left[\kappa \left(\frac{|H|^2}{n-1}-|A|^2\right)-\frac{n-1}{2}|A^-|^2\right]&=6|A^-|^2\left[\kappa\e|H|^2-|H|\lambda_1+\frac{n-1}{2}\lambda_1^2\right].
	\end{align*}
Normalise $|H|=1$. At $F_\e=0$, use \eqref{eqn_contact} in the reaction of \eqref{eqn_F}. Equations \eqref{rangeoflambda1}, \eqref{eqn_T-lower} and \eqref{eqn_norm-reactions} give
	\begin{align}\label{eqn_reaction-general}
	\mathcal R(F_\e)&\geq |A^-|^2\Bigg[\frac{n-2}{n}-(n-1)|\mathring h|^2+2\kappa\left(\frac4n-\frac3{n-1}\right)\\
	&-\left(4+\frac1n\right)\lambda_1+2(n-1)\lambda_1^2-\bigl(1-(n-1)\lambda_1\bigr)\sqrt{\frac{n-1}{n}}\,|\mathring h|+6\kappa\e\Bigg].\nonumber
	\end{align}
At $F_\e=0$, normalise $|H|=1$. By \eqref{rangeoflambda1} and \eqref{eqn_F-tracefree}, we have
	\begin{align*}
	0&<\lambda_1\le\frac1n, \qquad |\mathring h|^2\le\frac1{3n}.
	\end{align*}
Hence
	\begin{align*}
	\frac{n-2}{n}-(n-1)|\mathring h|^2&\ge\frac{n-2}{n}-\frac{n-1}{3n}=\frac23-\frac{5}{3n}.
	\end{align*}
For the chosen value $\kappa=14/5$,
	\begin{align*}
	2\kappa\left(\frac4n-\frac3{n-1}\right)&=\frac{28(n-4)}{5n(n-1)}.
	\end{align*}
Next, since $0<\lambda_1\le1/n$,
	\begin{align*}
	-\left(4+\frac1n\right)\lambda_1+2(n-1)\lambda_1^2&\ge-5\lambda_1\ge-\frac5n.
	\end{align*}
Also, $0<1-(n-1)\lambda_1\le1$, so
	\begin{align*}
	-\bigl(1-(n-1)\lambda_1\bigr)\sqrt{\frac{n-1}{n}}\,|\mathring h|&\ge-\sqrt{\frac{n-1}{n}}\,|\mathring h|\ge-\frac{\sqrt{n-1}}{\sqrt3\,n}.
	\end{align*}
Therefore,
	\begin{align*}
	\mathcal R(F_\e)&\geq \left(\frac23-\frac{16}{15n}-\frac{84}{5n(n-1)} -\frac{\sqrt{n-1}}{\sqrt3\,n}+6\kappa\e\right)|A^-|^2,
	\end{align*}
with
	\begin{align*}
	0&< c_1\le\frac23-\frac{16}{15n}-\frac{84}{5n(n-1)} -\frac{\sqrt{n-1}}{\sqrt3\,n},\quad n\ge10.
	\end{align*}
Restoring scale gives \eqref{finalestimateR}.
\end{proof}

\subsection{Gradient terms}
For the following algebraic estimates, we normalise $|H|=1$, so that $\operatorname{tr}h=1$. Let $T$ be a fully symmetric tensor of rank three on $ \mathbb{R}^n$ with trace $w_k=\sum_iT_{kii}$. Write $T=E(w)+T^0$, where the tensor $T^0$ is symmetric and trace free. The symmetric tensor of smallest norm with trace $w$ is
	\begin{align}\label{E(w)}
	E(w)_{kij}&=\frac{\delta_{ij}w_k+\delta_{ki}w_j+\delta_{kj}w_i}{n+2}.
	\end{align}
Then
	\begin{align}\label{eqn_codazzi-trace}
	|T|^2&=|E(w)|^2+|T^0|^2\ge\frac3{n+2}|w|^2,
	\end{align}
	\begin{align}\label{eqn_codazzi-projection}
	|T-h\otimes w|^2&\ge\frac23\sum_j\left(|h|^2-\lambda_j^2+\frac{2(1-\lambda_j)^2}{n+2}\right)w_j^2,
	\end{align}
where $(h\otimes w)_{kij}=h_{ij}w_k$.
First split $h \otimes w$ into its symmetric part $\operatorname{Sym}(h \otimes w)$ and its nonsymmetric part $h\otimes w-\operatorname{Sym}(h \otimes w)$. Since $T$ is symmetric, these two pieces are orthogonal, so
	\begin{align*}
	|T-h \otimes w|^2&=|T- \operatorname{Sym}(h \otimes w)|^2+|h \otimes w-\operatorname{Sym}(h \otimes w)|^2.
	\end{align*}
For the nonsymmetric part, a direct calculation gives
	\begin{align}\label{firstpart}
	|h \otimes w- \operatorname{Sym}(h \otimes w)|^2&=\frac{2}{3}\left(|h|^2|w|^2-|h w|^2\right)=\frac{2}{3} \sum_j\left(|h|^2-\lambda_j^2\right) w_j^2,
	\end{align}
where the last equality is assuming $h$ is diagonal with eigenvalues $\lambda_j$. For the symmetric part,
	\begin{align*}
	\operatorname{tr} \operatorname{Sym}(h \otimes w)_{kii}&= \frac{1}{3}\left(h_{ii}w_k+h_{ki}w_i+h_{ki}w_i \right)=\frac{w_k+2 h_{ki} w_i}{3},
	\end{align*}
because $\operatorname{tr} h=1$. Therefore its trace free part is
	\begin{align*}
	S^0&=\operatorname{Sym}(h \otimes w)-E\left(\frac{w+2 h w}{3}\right).
	\end{align*}
Since $T^0$ is arbitrary, to minimise the distance we choose $T^0$ equal to exactly this trace free part. What remains is only the difference of the trace parts
	\begin{align*}
	E(w)-E\left(\frac{w+2 h w}{3}\right)&=E\left(\frac{2(w-h w)}{3}\right).
	\end{align*}
Since the trace and trace free parts are orthogonal,
	\begin{align*}
	|T-\operatorname{Sym}(h \otimes w)|^2&=\left|E\left(\frac{2(w-h w)}{3}\right)\right|^2+\left|T^0-S^0\right|^2\geq\left|E\left(\frac{2(w-h w)}{3}\right)\right|^2.
	\end{align*}
But
	\begin{align}\label{secondpart}
	\left|E\left(\frac{2(w-h w)}{3}\right)\right|^2&=\frac{4}{3(n+2)}|w-h w|^2 =\frac{4}{3(n+2)} \sum_j\left(1-\lambda_j\right)^2 w_j^2,
	\end{align}
Adding \eqref{firstpart} and \eqref{secondpart} gives \eqref{eqn_codazzi-projection}.

At an interior zero of $F_\e$, with $\kappa=14/5$, the estimate \eqref{eqn_F-tracefree} and the normalisation $|H|=1$ give
	\begin{align*}
	|\mathring h|^2&\le\frac1{3n}.
	\end{align*}
Consequently,
	\begin{align*}
	\lambda_n-\frac1n&\le\sqrt{\frac{n-1}{n}}\,|\mathring h| \le\frac{\sqrt{(n-1)/3}}n.
	\end{align*}
For $n\ge10$, it follows that
	\begin{align*}
	1-\lambda_n&\ge\frac{n-1-\sqrt{(n-1)/3}}n \ge\frac{4(n-1)}{5n}.
	\end{align*}
Using
	\begin{align*}
	|h|^2-\lambda_j^2&=\sum_{i\ne j}\lambda_i^2 \ge\frac{(1-\lambda_j)^2}{n-1}
	\end{align*}
in \eqref{eqn_codazzi-projection}, we obtain
	\begin{align}\label{eqn_codazzi-final}
	|T-h\otimes w|^2&\ge\frac{2n}{(n-1)(n+2)} \sum_j(1-\lambda_j)^2w_j^2\ge\frac{2n(1-\lambda_n)^2}{(n-1)(n+2)}|w|^2\\
	&\ge\frac{32(n-1)}{25n(n+2)}|w|^2.\nonumber
	\end{align}
\begin{lemma}[Estimate for the gradient terms]
Let $n\ge10$, $\kappa=14/5$ and $\e>0$. At an interior zero of $F_\e$ with $|H|>0$, we have
	\begin{align}\label{finalestimateG}
	\mathcal G(F_\e)&\geq c_2|H|^2|\nabla^\perp\nu|^2,
	\end{align}
for some $c_2>0$.
\end{lemma}
\begin{proof}
Let $F_\e=0$. We normalise $|H|=1$, diagonalise $h$ and choose $e=e_1$. Also, denote $a=(n-1)/2$ and
	\begin{align*}
	(T_\alpha)_{kij}&=\ip{\nabla_k A_{ij}}{\nu_\alpha},\qquad p_k=\sum_i(T_1)_{kii}=\nabla_k |H|,\\
	(w_\alpha)_k&=\sum_i(T_\alpha)_{kii}=\ip{\nabla_k^\perp H}{\nu_\alpha},\qquad\alpha>1.
	\end{align*}
Each $T_\alpha$ is fully symmetric by Codazzi. Define
	\begin{align*}
	W_{kij}&=\sum_{\alpha>1}(A^-_\alpha)_{ij}(w_\alpha)_k, \qquad |w|^2=\sum_{\alpha>1}|w_\alpha|^2.
	\end{align*}
Then $(\nabla_k h)_{ij}=(T_1)_{kij}+W_{kij}$ at the point.
\begin{claim}
The gradient terms are given by
	\begin{align}\label{eqn_gradient-exact}
	\frac12\mathcal G(F_\e)&=\kappa|T_1|^2-2a\ip{T_1}{W}-a|W|^2+a|\nabla h_{11}|^2-p\cdot\nabla h_{11}\\
	&+p\cdot\left(2a\ip{h}{W}+(1-2a\lambda_1)W_{11}\right)-\kappa\left(\frac{1}{n-1}-\e\right)|p|^2\nonumber\\
	&+\sum_{\alpha>1}\Bigg\{\kappa|T_\alpha|^2+a|T_\alpha-h\otimes w_\alpha|^2-(1-2a\lambda_1)\sum_k(T_\alpha)_{k11}(w_\alpha)_k\nonumber\\
	&-\left(a\lambda_1^2+\kappa\left(\frac{1}{n-1}-\e\right)\right)|w_\alpha|^2 \Bigg\}\nonumber.
	\end{align}
\end{claim}
\begin{proof}
For a fixed spatial direction, consider the linear variation $A_\alpha(s)=A_\alpha+sT_\alpha$, regarding $T_\alpha$ as a matrix with its spatial derivative index fixed. Therefore, $H(s)=H+s \operatorname{tr} T$. Decompose the normal vector $\operatorname{tr} T$ in the normal frame $\left\{\nu_\alpha\right\}$,
	\begin{align*}
	\operatorname{tr} T&=p \nu_1+\sum_{\alpha>1} w_\alpha \nu_\alpha.
	\end{align*}
Therefore
	\begin{align*}
	H(s)&=\nu_1+s\left(p \nu_1+\sum_{\alpha>1} w_\alpha \nu_\alpha\right)=(1+s p) \nu_1+s \sum_{\alpha>1} w_\alpha \nu_\alpha.
	\end{align*}
At $s=0$ we have
	\begin{align*}
	\frac{d}{ds}|H|&=p,\\
	\frac{d^2}{ds^2}|H|&=\sum_{\alpha>1}w_\alpha^2,\\
	\frac{d}{ds}h&= \frac{d}{ds}\langle A(s),\nu(s)\rangle=T_1+\sum_{\alpha>1}w_\alpha A^-_\alpha,\\
	\frac{d^2}{ds^2}h&=2\sum_{\alpha>1}w_\alpha T_\alpha-2p\sum_{\alpha>1}w_\alpha A^-_\alpha-h\sum_{\alpha>1}w_\alpha^2.
	\end{align*}
Therefore, by applying $-\tfrac12d^2/ds^2$ to
	\begin{align*}
	|H|h_{11}-a h_{11}^2+\kappa\left(\frac{1}{n-1}-\e\right)|H|^2-(a+\kappa)|A|^2+a|h|^2
	\end{align*}
which is \eqref{eqn_F} rearranged and sum over spatial directions, this yields \eqref{eqn_gradient-exact}, including all derivatives of $H/|H|$.
\end{proof}
We begin by estimating the terms in \eqref{eqn_gradient-exact}. Note that
	\begin{align}\label{psquared}
	a|\nabla h_{11}|^2-p\cdot\nabla h_{11}&=a\left|\nabla h_{11}-\frac{p}{2a}\right|^2-\frac{|p|^2}{4a}\ge-\frac{|p|^2}{4a}.
	\end{align}
Every symmetric trace free matrix $B$ satisfies $\|B\|_{\mathrm{op}}\le\sqrt{(n-1)/n}\,|B|$. Consequently,
	\begin{align}\label{Aw}
	\left|\sum_{\alpha>1}A^-_\alpha w_\alpha\right|&\leq\left(\sum_\alpha\left\|A_\alpha^-\right\|_{\text{ op }}^2\right)^{1 / 2}\left(\sum_\alpha\left|w_\alpha\right|^2\right)^{1 / 2}\leq \sqrt{\frac{n-1}{n}}\left|A^-\right||w|\\
	&\le\sqrt{\frac{n-1}{n}}|A^-|\,|w|.\nonumber
	\end{align}
Similarly, using $\left|\left(A_\alpha^-\right)_{11}\right| \leq\sqrt{(n-1)/n}\left|A_\alpha^-\right|$ and Cauchy-Schwarz, we obtain
	\begin{align*}
	\left|W_{11}\right|&=\left|\sum_{\alpha>1}\left(A_\alpha^-\right)_{11} w_\alpha\right| \leq \left(\sum_{\alpha>1}\left|\left(A_\alpha^-\right)_{11}\right|^2\right)^{1 / 2}\left(\sum_{\alpha>1}\left|w_\alpha\right|^2\right)^{1 / 2}\leq \sqrt{\frac{n-1}{n}}\left|A^-\right||w|.
	\end{align*}
Lastly, since $A_\alpha^-$ is trace free, we have $\left\langle h, A_\alpha^-\right\rangle=\left\langle\mathring{h}, A_\alpha^-\right\rangle$. Hence $|\langle h, W\rangle| \leq|\mathring{h}| |A^-||w|$ and therefore
	\begin{align}\label{hWterm}
	2 a p\langle h, W\rangle&\geq-(n-1)|\mathring{h}| |A^-||p||w|.
	\end{align}
For each $\alpha>1$, set $V=T_\alpha-h\otimes w_\alpha$. The term inside the sum in \eqref{eqn_gradient-exact} is
	\begin{align*}
	\kappa|T_\alpha|^2+a|V|^2 -(1-2a\lambda_1)\ip{V_{11}}{w_\alpha}-\left(\kappa\left(\frac{1}{n-1}-\e\right)+\lambda_1-a\lambda_1^2\right)|w_\alpha|^2.
	\end{align*}
Since $0<\lambda_1\le1/n$, we have $0<1-2a\lambda_1\le1$ and $\lambda_1-a\lambda_1^2\le1/n$. Young's inequality gives
	\begin{align*}
	-(1-2a\lambda_1)\ip{V_{11}}{w_\alpha}&\ge-\frac32|V|^2-\frac16|w_\alpha|^2.
	\end{align*}
Applying \eqref{eqn_codazzi-trace} to $T_\alpha$ and \eqref{eqn_codazzi-final} to $V$, the above is bounded below by
	\begin{align*}
	\left[ \frac{14}{5}\left(\frac3{n+2}-\frac1{n-1}\right)+\frac{16(n-4)(n-1)}{25n(n+2)}-\frac1n-\frac16 \right]|w_\alpha|^2&\ge\frac25|w_\alpha|^2.
	\end{align*}
Subtracting $\frac25|w_\alpha|^2$ from both sides, we get
	\begin{align*}
	\frac{11n^3+29n^2-1216n-84}{150n(n-1)(n+2)}|w_\alpha|^2&\ge0.
	\end{align*}
The numerator is increasing for $n\ge10$ and is positive at $n=10$. For the remaining terms in \eqref{eqn_gradient-exact}, we use \eqref{psquared}, \eqref{Aw}, and \eqref{hWterm}, together with $|W|^2 \leq\left|A^-\right|^2|w|^2$. We also use the slightly coarser estimate
	\begin{align*}
	\kappa\left|T_1\right|^2-2 a\left\langle T_1, W\right\rangle&\geq \frac{2 \kappa}{3}\left|T_1\right|^2-\frac{3 a^2}{\kappa}|W|^2,
	\end{align*}
which follows directly from Young's inequality. Applying \eqref{eqn_codazzi-trace} to $T_1$, and using $2 a=n-1$ and $\kappa=14 / 5$, we obtain
	\begin{align}\label{eqn_gradient-bound}
	\frac12\mathcal G(F_\e)&\geq \left(\frac{28}{5(n+2)}-\frac{33}{10(n-1)}\right)|p|^2 +\left(\frac25-\frac{(n-1)(15n+13)}{56}|A^-|^2\right)|w|^2\\
	&-\left((n-1)|\mathring h|+\sqrt{\frac{n-1}{n}}\right)|A^-|\,|p|\,|w|.\nonumber
	\end{align}
For $n\ge10$,
	\begin{align*}
	\frac{28}{5(n+2)}-\frac{33}{10(n-1)}-\frac1n&=\frac{(n-10)(13n-2)}{10n(n-1)(n+2)}\ge0.
	\end{align*}
Equations \eqref{eqn_contact-smallness}, \eqref{eqn_contact-product} and \eqref{eqn_contact-uniform} give
	\begin{align*}
	\frac{(n-1)(15n+13)}{56}|A^-|^2&\le\frac{17(15n+13)}{896n}:=C_1,\\
	\frac n4\left((n-1)|\mathring h|+\sqrt{\frac{n-1}{n}}\right)^2|A^-|^2&\le\frac{17}{32n}+\frac{289(5n+23)}{57344n}:=C_2,
	\end{align*}
with $C_1+C_2<\frac{2}{5}$. Both $C_1$ and $C_2$ are decreasing functions of $n$, so it suffices to check the inequality at $n=10$. For the second estimate use $(u+v)^2\le2u^2+2v^2$ and the product bound in \eqref{eqn_contact-product}. Completing the square in $|p|$ in \eqref{eqn_gradient-bound} we have
	\begin{align*}
	\mathcal G(F_\e)&\ge2\left(\frac25-C_1-C_2\right)|w|^2=c_2|w|^2.
	\end{align*}
The scalar $F_{\e,e}$ is homogeneous of degree two in $A$, so its Hessian in $A$ is homogeneous of degree zero.
Restoring the original scale and using $|w|^2=|H|^2|\nabla^\perp\nu|^2$ gives \eqref{finalestimateG}.
\end{proof}

\subsection{Completion of proof of Lemma \ref{lem_interior}}
\begin{proof}
For $n\ge10$, \eqref{finalestimateR} and \eqref{finalestimateG} give the required reaction and gradient bounds. For $4\leq n\le9$, we retain more of the gradient contribution in \eqref{eqn_gradient-exact}. Normalise $|H|=1$ and put $a=(n-1)/2$. Split $T_1=E(p)+T_1^0$ using \eqref{E(w)}, so that
	\begin{align*}
	\langle E(p),W\rangle&=\frac{2}{n+2}p\cdot\sum_{\alpha>1}A^-_\alpha w_\alpha.
	\end{align*}
Complete the squares in $T_1^0$, $\nabla h_{11}$ and $(T_\alpha-h\otimes w_\alpha)_{11}$, using $\lambda_1-a\lambda_1^2+(1-2a\lambda_1)^2/(4a)=1/(4a)$ in the last square. With $P_n=(51n-150)/(10(n-1)(n+2))$, completing the square in $|p|$ and using \eqref{Aw} and \eqref{hWterm} gives
	\begin{align*}
	\frac12\mathcal G(F_\e)&\ge\left[P_n-\left(a+\frac{a^2}{\kappa} +\frac{\bigl((n-1)|\mathring h|+3\sqrt{n(n-1)}/(n+2)\bigr)^2}{4P_n}\right)|A^-|^2\right]|w|^2.
	\end{align*}
For $5\leq n\le9$, apply $(s+t)^2\le5s^2+5t^2/4$ and \eqref{eqn_contact-smallness}. The quantity subtracted from $P_n$ is at most $(b_n+d_n x)(U_n-\beta_n x)$, where
	\begin{align*}
	x&=|\mathring h|^2,\qquad b_n=a+\frac{a^2}{\kappa}+\frac{45n(n-1)}{16P_n(n+2)^2}, \qquad d_n=\frac{5(n-1)^2}{4P_n}.
	\end{align*}
For each of these five dimensions $d_nU_n\leq b_n\beta_n$ and $P_n-b_nU_n\ge23/980$, so the product is decreasing for $x\ge0$ and $\mathcal G(F_\e)\ge(23/490)|w|^2$. For the reaction terms, use $|\mathring h|^2\le1/(3n)$ and $\sqrt{(n-1)/n}\,|\mathring h|\le1/4$ in \eqref{eqn_reaction-general}. Completing the square in $\lambda_1$ gives
	\begin{align*}
	\mathcal R(F_\e)&\geq |A^-|^2\left[\frac{2n-5}{3n}+\frac{28(n-4)}{5n(n-1)}-\frac14 -\frac{\bigl(4+1/n-(n-1)/4\bigr)^2}{8(n-1)}\right] \geq \frac{13}{300}|A^-|^2.
	\end{align*}

For $n=4$, the trace decomposition of $F_0$ gives $\lambda_1^2\le2\kappa\e$ at the zero. In an eigenframe, $|A|^2\le1/3$ gives
	\begin{align*}
	\left|h-\operatorname{diag}(0,1/3,1/3,1/3)\right|^2+|A^-|^2&\le\frac{2\lambda_1}{3}.
	\end{align*}
Thus the gradient bound above gives $\mathcal G(F_\e)\geq (1/2)|w|^2$ for sufficiently small $\e$, since $P_4=3/10$. For the reaction terms, substitute \eqref{eqn_contact} into the reaction of \eqref{eqn_F} using \eqref{normevolutions}. The mixed Gram and commutator terms have coefficients $-4\kappa$ and $-(2a+4\kappa)$ and depend, respectively, on the diagonal entries of $A^-_\alpha$ and on the entries off its diagonal. Using \eqref{eqn_T-lower} for $(1-3\lambda_1)\mathcal T_{11}$ and bounding the purely normal quartic terms as in \eqref{eqn_norm-reactions}, we obtain
	\begin{align*}
	\mathcal R(F_\e)\geq |A^-|^2\Big[&-\frac12+2(a+\kappa)|h|^2-6(a+\kappa)|A^-|^2\\
	&-\max\!\left\{4\kappa|\mathring h|^2,(2a+4\kappa)(\lambda_4-\lambda_1)^2\right\}\Big].
	\end{align*}
The preceding closeness estimate shows that the bracket tends uniformly to $71/90$ as $\e\to 0$. This proves the reaction bound for sufficiently small $\e$.

Restoring scale in the bounds in lower dimensions gives the terms on the right of \eqref{finalestimate}. The local extension of $e$ contributes no additional terms. At the point, $\nabla e=0$, $h(e,\cdot)=\lambda_1g(e,\cdot)$ and $\langle\Delta e,e\rangle=-|\nabla e|^2=0$.
\end{proof}

The interior maximum principle argument applies on any time interval where $|H|>0$ and the hypotheses of Lemma \ref{lem_interior} hold at a possible zero. Starting with $F_\e>0$, replace $\e$ by $\e-\sigma t$, with $\sigma>0$ small enough that this coefficient stays positive on the interval under consideration. At a first interior zero of the perturbed quantity, Lemma \ref{lem_modified-principal-convexity} gives a smooth upper support with a local minimum equal to zero. Lemma \ref{lem_interior}, applied with the instantaneous value of the coefficient, gives the additional time derivative $\kappa\sigma|H|^2>0$, contradicting the parabolic minimum principle. Letting $\sigma\to 0$ gives the interior maximum principle argument for $F_\e\ge0$.

\subsection{Preservation and convergence for closed flows}

The new preservation result supplies the convexity estimate needed to apply the existing derivative estimates and blowup classification of Lynch--Nguyen. Throughout this subsection, $n\ge4$ and $\kappa=14/5$ as in \eqref{eqn_kappa-choice}.

\begin{proposition}[Preservation and derivative estimates for closed flows]
\label{prop_closed-preservation-gradient}
Let $F:M^n\times[0,T)\to\mathbb R^{n+m}$ be a smooth MCF of a compact connected manifold without boundary. Fix $\e>0$, sufficiently small if $n=4$. Suppose $F_\e>0$ initially. If $4\leq n\le7$, assume in addition that $|A|^2<\frac{4}{3n}|H|^2$ on $M_0$. Then $F_\e>0$ throughout the smooth flow. Moreover, there are constants $a>0$, $\frac{1}{n}<c<\frac{4}{3n}$, $\eta>0$ and $C<\infty$ such that
	\begin{align}\label{eq_closed-curvature-bounds}
	h&\ge\eta|H|g, \qquad |A|^2+a\leq c|H|^2,
	\end{align}
and
	\begin{align}\label{eq_closed-derivative-estimates}
	|\nabla A|^2&\leq C|H|^4, \qquad |\nabla^2A|^2\leq C|H|^6.
	\end{align}
\end{proposition}

\begin{proof}
By the convention following \eqref{eqn_F}, $F_\e>0$ implies that $H$ is nowhere zero on $M_0$. The strict bound $|A|^2<\frac{4}{3n}|H|^2$ holds by hypothesis for $4\leq n\le7$ and by \eqref{eq_modified-strict-coarse} for $n\ge8$. Choose $c$ strictly between $\max_{M_0}|A|^2/|H|^2$ and $\frac{4}{3n}$. The continuous function $c|H|^2-|A|^2$ is strictly positive on $M_0$, so compactness gives
	\begin{align*}
	a&=\min_{M_0}\bigl(c|H|^2-|A|^2\bigr)>0.
	\end{align*}
Lynch--Nguyen \cite[Lemma 2.2]{LNSurgery2025} show that the flow preserves $|A|^2+a\leq c|H|^2$. Since $|A|^2\ge0$, it follows that
	\begin{align*}
	|H|^2&\ge\frac ac>0
	\end{align*}
throughout the smooth flow. This supplies the nonvanishing of $H$ and the quadratic pinching needed in Lemma \ref{lem_interior}.

By compactness, $\min_{M_0}F_\e/|H|^2>0$. Choose $\e'>\e$ sufficiently close to $\e$ that $F_{\e'}=F_\e-\kappa(\e'-\e)|H|^2>0$ initially and, in dimension four, $\e'$ is still within the range of Lemma \ref{lem_interior}. Since $M$ has no boundary, the maximum principle argument above gives $F_{\e'}\ge0$ throughout the flow. Thus $F_\e\ge\kappa(\e'-\e)|H|^2>0$. By Lemma \ref{lem_modified-principal-convexity}, we also have
	\begin{align*}
	h&\ge\eta|H|g, \qquad \eta=\frac{\kappa\e}{1+2\kappa/(n-1)}>0.
	\end{align*}
This proves \eqref{eq_closed-curvature-bounds}.

Lynch--Nguyen's gradient estimate and smooth flow second derivative argument \cite[Theorem 4.1 and proof of Theorem 4.3]{LNSurgery2025} apply under the preserved quadratic pinching with $c<\frac{4}{3n}$. The latter argument uses only the quadratic gap and the gradient estimate, so it also applies in dimension four. Smoothness on a compact initial time interval extends the estimates to the initial time after enlarging the constant. Thus
	\begin{align*}
	|\nabla A|^2&\leq C(|A|^4+1), \qquad |\nabla^2A|^2\leq C(|A|^6+1).
	\end{align*}
Using $|A|^2\leq c|H|^2$ and the positive lower bound for $|H|$, we absorb the constant terms and obtain \eqref{eq_closed-derivative-estimates}.
\end{proof}

\begin{proof}[Proof of Theorem \ref{thm_closed-convergence}]
Proposition \ref{prop_closed-preservation-gradient} gives preservation of $F_\e>0$ and the curvature and derivative estimates. For $n=4$, the convergence follows directly from \cite[Theorem 1]{AndrewsBaker}. For $n\ge5$, let $c_n$ be the pinching coefficient in \cite[Theorem 1.2]{LN}, equal to $\frac{4}{3n}$ for $n\ge8$ and to the stated coefficient for $5\leq n\le7$. In the proof of Proposition \ref{prop_closed-preservation-gradient}, we can therefore choose $c<c_n$.

\smallskip\noindent\emph{Singularity formation in finite time and rescaling.} The identity
	\begin{align*}
	(\partial_t-\Delta)|F|^2&=-2n
	\end{align*}
and the maximum principle imply
	\begin{align*}
	\max_M|F(\cdot,t)|^2&\le\max_M|F(\cdot,0)|^2-2nt.
	\end{align*}
Thus $T<\infty$, and the extension criterion gives unbounded curvature as $t\to T$, see the continuation argument in \cite[Section 4]{AndrewsBaker}. By \eqref{eq_closed-curvature-bounds}, $|H|$ is also unbounded.

We follow the rescaling and compactness argument in the proof of \cite[Corollary 5.1]{LN}. Choose $t_j\to T$ and $p_j\in M$ such that
	\begin{align*}
	|H|(p_j,t_j)&=\max_{M\times[0,t_j]}|H| \longrightarrow\infty.
	\end{align*}
Set $r_j=|H|(p_j,t_j)^{-1}$ and define
	\begin{align*}
	F_j(p,\tau)&=r_j^{-1}\bigl(F(p,t_j+r_j^2\tau)-F(p_j,t_j)\bigr), \qquad -t_j/r_j^2\le\tau\le0.
	\end{align*}
The rescaled flows satisfy
	\begin{align*}
	|H_j|&\le1,\qquad |A_j|^2\leq c, \qquad |H_j|(p_j,0)=1.
	\end{align*}
The estimates \eqref{eq_closed-derivative-estimates}, the neighbourhood estimates at the curvature scale in \cite[Corollary 4.4 and Lemma 4.5]{LNSurgery2025}, and standard interior estimates give bounds for all derivatives on compact parabolic neighbourhoods. Smooth pointed compactness, as used in \cite[Section 5]{LN}, therefore gives a subsequence converging to a complete connected ancient flow
	\begin{align*}
	F_\infty:M_\infty\times(-\infty,0]&\longrightarrow\mathbb R^{n+m},
	\end{align*}
with $|H_\infty|(p_\infty,0)=1$. In particular, the limit is nonflat. Rescaling at times when the running maximum is attained requires no Type I assumption; we use only the rescaling and compactness argument from the proof of \cite[Corollary 5.1]{LN}.

\smallskip\noindent\emph{Compactness and roundness of the ancient limit.} Since $c<c_n$, the bound \eqref{eq_closed-curvature-bounds} satisfies the quadratic pinching hypothesis of \cite[Theorem 1.2]{LN}. That theorem, incorporating Naff's planarity estimate \cite{Naff2022}, shows that the limiting flow consists of convex hypersurfaces bounding convex sets in a fixed $(n+1)$-dimensional affine subspace.

The preserved principal convexity estimate passes to the limit:
	\begin{align*}
	h_\infty&\ge\eta|H_\infty|g_\infty.
	\end{align*}
Nonflatness and the strong maximum principle imply strict convexity. Completeness and strict convexity identify the limiting immersion at each fixed time with the boundary of a convex set. Hamilton's compactness theorem \cite[Theorem 1.1]{HamiltonCompactness1994} then implies that the submanifold is compact at every time. Extend the limit smoothly forward from $\tau=0$ by Huisken's convex hypersurface flow \cite{Huisken1984}. The tensor maximum principle preserves the same convexity bound up to the finite extinction time $T_\infty>0$. With the shifted time $s=\tau-T_\infty$, the singular time is $s=0$, so Huisken-Sinestrari \cite[Theorem 3.2]{HuiskenSinestrari2015} imply that the limiting flow is a family of shrinking round spheres. In particular,
	\begin{align*}
	\frac{|A_\infty|^2}{|H_\infty|^2}&=\frac1n.
	\end{align*}

\smallskip\noindent\emph{Global convergence and critical quadratic pinching at a late time.} It remains to pass from the compact pointed limit to the whole original submanifold. At $\tau=0$, the convergence maps from the pointed compactness construction are defined on an exhaustion of $M_\infty$. Since $M_\infty$ is compact, they are defined on all of $M_\infty$ for sufficiently large $j$. Each is a diffeomorphism onto an open subset of $M$; its image is also closed, by compactness. Connectedness of $M$ implies that the image is all of $M$. Thus the convergence at $\tau=0$ is global after reparametrisation.

Since the limiting sphere has positive mean curvature, global smooth convergence gives
	\begin{align*}
	\sup_M\frac{|A|^2}{|H|^2}(\cdot,t_j)&\longrightarrow\frac1n.
	\end{align*}
Choose $c_*$ with $\frac{1}{n}<c_*<\frac{1}{n-1}$. For sufficiently large $j$, the entire submanifold $M_{t_j}$ satisfies
	\begin{align*}
	|A|^2&\leq c_*|H|^2.
	\end{align*}
Restarting the flow at $t_j$, the Andrews--Baker convergence theorem \cite[Theorem 1]{AndrewsBaker} gives convergence of the original flow to a round point, with smooth convergence under the normalisation in the statement. This completes the proof.
\end{proof}

\begin{corollary}
A compact connected submanifold without boundary satisfying $|A|^2<\frac{1}{n-1}|H|^2$ initially evolves to a round point under MCF when $n\ge4$. The same conclusion holds under nonstrict critical quadratic pinching if $h$ is initially defined and positive definite everywhere.
\end{corollary}
\begin{proof}
Lemma \ref{lem_initial-positivity} gives $F_\e>0$ for sufficiently small $\e>0$ with the fixed choice of $\kappa$. For $5\leq n\le7$, we have $\frac{1}{n-1}<c_n$, so Theorem \ref{thm_closed-convergence} applies in every dimension $n\ge5$. In dimension four the convergence follows from \cite[Theorem 1]{AndrewsBaker}, which allows the nonstrict quadratic inequality when $H$ is nowhere zero.
\end{proof}

\section{Quadratic Pinching}\label{sec:quadratic-pinching}

The spherical boundary imposes a lower bound on the coefficient for which the maximum principle preserves ordinary quadratic pinching. The conormal identity below gives the range in Proposition \ref{prop_fb-quadratic-preservation} and explains the boundary restriction that motivates quadratic convexity. Throughout this section, $M_t$ is a smooth compact MCF of dimension $n\ge2$ in $\overline B_R^{n+m}$, with orthogonal free boundary on $\mathbb S_R^{n+m-1}$. We use the boundary frame of Section \ref{sec:notation}, with $e_n=\mu=F/R$.

\begin{lemma}[Spherical boundary identities]\label{nablaH}
Along the free boundary,
	\begin{align*}
	A_{a\mu}&=0,\qquad R\nabla_{e_a}\mu=e_a,\qquad R\nabla^\perp_\mu H=H,
	\end{align*}
	\begin{align}\label{eqn_nablamuAab}
	R(\nabla_\mu A)_{ab}&=\delta_{ab}A_{\mu\mu}-A_{ab},
	\end{align}
	\begin{align}\label{eqn_nablamuAmumu}
	R(\nabla_\mu A)_{\mu\mu}&=2H-nA_{\mu\mu}.
	\end{align}
\end{lemma}

\begin{proof}
Differentiating $\mu=F/R$ along $e_a$ gives $D_{e_a}\mu=e_a/R$, whose tangent and normal components give the first two identities. Differentiating $\langle e_b,\mu\rangle=0$ then gives $\langle\nabla_{e_a}e_b,\mu\rangle=-\delta_{ab}/R$.

For the time derivative, write $\mu=dF(u)$ with $u$ a tangent field that varies with time along the boundary of $M$. Since $\partial_tF=H$, the normal component of $D_t(dF(u))$ is $\nabla^\perp_\mu H$. On the other hand, $\mu=F/R$ gives $D_t\mu=H/R$.

By Codazzi and the vanishing of $A_{\mu b}$ along the boundary,
	\begin{align*}
	(\nabla_\mu A)_{ab}&=(\nabla_a A)_{\mu b}=-A(\nabla_{e_a}\mu,e_b)-A(\mu,\nabla_{e_a}e_b)=\frac1R\bigl(\delta_{ab}A_{\mu\mu}-A_{ab}\bigr).
	\end{align*}
Taking the tangential trace gives $R\sum_a(\nabla_\mu A)_{aa}=nA_{\mu\mu}-H$. Subtracting this from $R\nabla^\perp_\mu H=H$ gives the last identity.
\end{proof}

Set $P=|A|^2-c|H|^2$. Contracting the boundary identities gives
	\begin{align*}
	\frac R2\,\partial_\mu|A|^2&=3\langle H,A_{\mu\mu}\rangle-n|A_{\mu\mu}|^2-|A|^2,\\
	\frac R2\,\partial_\mu|H|^2&=|H|^2.
	\end{align*}
Subtracting $c$ times the second identity and completing the square yields
	\begin{align}\label{eq_quadratic-boundary}
	\frac R2\,\partial_\mu P&=-P-n\left|A_{\mu\mu}-\frac{3}{2n}H\right|^2 +\left(\frac9{4n}-2c\right)|H|^2.
	\end{align}
For $0<c\leq\frac{4}{3n}$, the reaction estimate \eqref{eqn_norm-reactions2}, written with coefficient $c$, and the Kato inequality give, wherever $|H|>0$,
	\begin{align}\label{eq_fb-quadratic-evolution}
	(\partial_t-\Delta)P&\leq 2\bigl(|h|^2+3|A^-|^2\bigr)P +2\left(3c-\frac4n\right)|H|^2|A^-|^2.
	\end{align}
Here the gradient contribution is nonpositive because $c\leq\frac{4}{3n}<\frac{3}{n+2}$.

\begin{proposition}[Preservation of quadratic pinching]\label{prop_fb-quadratic-preservation}
If $P\leq0$ initially and $\frac{9}{8n}\leq c\leq\frac{4}{3n}$, then $P\leq0$ throughout the flow.
\end{proposition}

\begin{proof}
Where $P>0$, \eqref{eq_fb-quadratic-evolution} gives $(\partial_t-\Delta)P\leq6|A|^2P$. This bound extends across $H=0$ by continuity of the reaction terms and the same Kato inequality. The boundary identity \eqref{eq_quadratic-boundary} gives $\partial_\mu P\leq-2P/R$.

Fix $S<T$ and $C>6\sup_{M\times[0,S]}|A|^2$. For $\eta>0$, the function $P-\eta e^{Ct}$ is initially negative. At a first zero, $P=\eta e^{Ct}>0$. At an interior contact the evolution inequality makes $(\partial_t-\Delta)(P-\eta e^{Ct})<0$. At a boundary contact the outward derivative is strictly negative. Both contradict the differential inequalities at a first maximum. Letting $\eta\to0$ proves preservation on $[0,S]$, and $S<T$ was arbitrary. The argument includes both endpoints. In particular, $c=\frac{1}{n-1}$ lies in this range for $4\leq n\leq9$.
\end{proof}

\begin{lemma}[Additive quadratic gap]\label{lem_fb-additive-gap}
If $|A|^2<\frac{4}{3n}|H|^2$ initially, there are constants $\frac{5}{4n}<c_0<\frac{4}{3n}$ and $a>0$, such that
	\begin{align}\label{eq_fb-additive-gap}
	|A|^2+a&\leq c_0|H|^2,\qquad |H|^2\geq\frac{a}{c_0-\frac{1}{n}}>0
	\end{align}
throughout the flow.
\end{lemma}

\begin{proof}
The strict initial inequality gives $|H|>0$ on $M_0$. By compactness we can choose
	\begin{align*}
	\max\left\{\frac5{4n},\ \max_{M_0}\frac{|A|^2}{|H|^2}\right\}&<c_0<\frac4{3n}.
	\end{align*}
Set $P=|A|^2-c_0|H|^2$ and choose $0<a<\min_{M_0}(-P)$. At a possible first zero of $P+a$, the trace inequality gives $a\leq\left(c_0-\frac{1}{n}\right)|H|^2$, so $|H|>0$ there.

At an interior contact, \eqref{eq_fb-quadratic-evolution} with $P=-a$ gives
	\begin{align*}
	(\partial_t-\Delta)(P+a)&\leq-2a|h|^2<0.
	\end{align*}
At a boundary contact, \eqref{eq_quadratic-boundary} and the trace inequality give
	\begin{align*}
	\frac R2\,\partial_\mu(P+a)&\leq a+\left(\frac9{4n}-2c_0\right)|H|^2 \leq\left(\frac5{4n}-c_0\right)|H|^2<0.
	\end{align*}
The two signs rule out first contact. Thus $|A|^2+a\leq c_0|H|^2$ throughout the flow. Applying $|A|^2\geq \frac{1}{n}|H|^2$ once more gives the lower bound for $|H|^2$.
\end{proof}

\section{Planarity}\label{sec:planarity}

We extend the interior planarity estimate to spherical free boundary flow by controlling the conormal derivative of the planarity quotient. This boundary calculation gives the decay in the other normal directions used to identify singularity models and prove uniqueness of the limiting plane in Section \ref{sec:fb-convergence}.
Let $F:M^n\times[0,T)\to\overline B_R^{n+m}$, with $n\ge5$ and $m\ge2$, be a smooth compact MCF with spherical free boundary, satisfying
	\begin{align*}
	|A|^2&<\frac4{3n}|H|^2
	\end{align*}
initially. We prove that, for every fixed $0<\sigma\le1/30$,
	\begin{align*}
	|A^-|^2&\leq C|H|^{2-2\sigma}
	\end{align*}
throughout the flow, where $C$ depends on the initial submanifold, $R$ and $\sigma$.

By Lemma \ref{lem_fb-additive-gap}, there exist $\frac{5}{4n}<c_0<\frac{4}{3n}$ and $a>0$ such that
	\begin{align*}
	|A|^2+a&\leq c_0|H|^2.
	\end{align*}
In particular, $|H|>0$. Choose $c_0<c_1<\frac{4}{3n}$ and set
	\begin{align}
	Q&=c_1|H|^2-|A|^2 \geq a+(c_1-c_0)|H|^2>0,\qquad	P_\sigma=\frac{|A^-|^2}{Q^{1-\sigma}}\label{Psigma}.
	\end{align}
The local estimates in Lynch's proof of Theorem 1.1, specifically \cite[Propositions 5.1-5.2]{Lyn}, give
	\begin{align*}
	(\partial_t-\Delta)P_\sigma&\leq \frac{2(1-\sigma)}{Q} \langle\nabla Q,\nabla P_\sigma\rangle
	\end{align*}
in the interior. It remains to exclude a positive boundary maximum.

The quadratic boundary identity \eqref{eq_quadratic-boundary} gives
	\begin{align}\label{evolQ}
	\frac R2\partial_\mu Q&=-Q+n\left|A_{\mu\mu}-\frac{3}{2n}H\right|^2 +\left(2c_1-\frac9{4n}\right)|H|^2
	\end{align}
	\begin{align}\label{ineqQ}
	\frac R2\partial_\mu Q&\geq -Q,
	\end{align}
since $c_1>\frac{9}{8n}$. Moreover, $\nabla_\mu^\perp H=H/R$ from Lemma \ref{nablaH} implies
	\begin{align}\label{norotation}
	\nabla_\mu^\perp\nu&=0.
	\end{align}
Projecting the spherical boundary identities orthogonally to $\nu$, we obtain
	\begin{align}
	A^-_{a\mu}&=0,\\
	R(\nabla_\mu A^-)_{ab}&=\delta_{ab}A^-_{\mu\mu}-A^-_{ab}\label{nablaA-ab},\\
	R(\nabla_\mu A^-)_{\mu\mu}&=-nA^-_{\mu\mu}\label{nablaA-mm}.
	\end{align}
Since $A^-$ is trace free, $\sum_aA^-_{aa}=-A^-_{\mu\mu}$, and hence
	\begin{align}\label{partialmuA}
	\frac R2\partial_\mu|A^-|^2&=\sum_{a,b} \left\langle A^-_{ab}, \delta_{ab}A^-_{\mu\mu}-A^-_{ab}\right\rangle -n|A^-_{\mu\mu}|^2\\
	&=-|A^-|^2-n|A^-_{\mu\mu}|^2.\nonumber
	\end{align}
Combining these estimates gives
	\begin{align}\label{partialPsigma}
	\frac R2\partial_\mu P_\sigma&=\frac{R\,\partial_\mu|A^-|^2}{2Q^{1-\sigma}} -(1-\sigma)\frac{P_\sigma}{Q}\frac R2\partial_\mu Q\\
	&\le-\sigma P_\sigma -\frac{n|A^-_{\mu\mu}|^2}{Q^{1-\sigma}}.\nonumber
	\end{align}
Thus $\partial_\mu P_\sigma<0$ wherever $P_\sigma>0$, excluding a positive boundary maximum. The maximum principle therefore yields
	\begin{align*}
	\max_{M_t}P_\sigma&\le\max_{M_0}P_\sigma,
	\end{align*}
and
	\begin{align}\label{fbestimate}
	|A^-|^2&\leq C Q^{1-\sigma}\leq C|H|^{2-2\sigma}.
	\end{align}

Under a curvature rescaling with factors $\lambda_j\to\infty$, this becomes
	\begin{align*}
	|A_j^-|^2&\leq C\lambda_j^{-2\sigma}|H_j|^{2-2\sigma}.
	\end{align*}
Therefore, $A^-_\infty=0$ wherever $|H_\infty|>0$ in a smooth blowup limit. The strong maximum principle argument in \cite[proof of Corollary 2.6]{Naff2022}, applied to the limiting $Q$, shows that each connected limit is either flat or has $|H_\infty|>0$ everywhere. For boundary limits, apply this argument after reflection across the limiting support hyperplane, as described in Section \ref{sec:fb-convergence}.

In the nonflat case, the inherited pinching gives
	\begin{align*}
	h_\infty&\leq \sqrt{c_0}\,|H_\infty|g_\infty <|H_\infty|g_\infty.
	\end{align*}
Naff's Proposition 2.5 \cite{Naff2022} therefore places the limit in an affine $(n+1)$-plane. Thus every connected smooth curvature blowup is codimension one.

\section{Gradient Estimates}\label{sec:gradient-estimates}
We construct explicit correction tensors $B$ and $\widetilde B$ that give the conormal inequalities needed for the first two derivative estimates in arbitrary codimension. This construction extends the curvature estimates up to the spherical free boundary. Both corrections are defined algebraically on the whole evolving submanifold. We first work with $R=1$ and use the time connection, which preserves the metric. On tensors, we write $\mathcal L=\nabla_{\partial_t}-\Delta$. On scalar functions this is $\partial_t-\Delta$. Constants denoted by $C$ may change between lines.

\subsection{The first derivative correction}

Set $v=F^\top$, $w=F^\perp$ and in an orthonormal tangent frame, define
	\begin{align*}
	v_i&=\langle F,e_i\rangle,\qquad \pi_{ij}=\delta_{ij}-v_iv_j,\qquad A_{vv}=\sum_{p,q}v_pv_qA_{pq}.
	\end{align*}
The tensor $\pi$ is the orthogonal projection onto $T\partial M_t$ at the boundary. In the interior it is simply the coefficient defined by this formula. The Gauss and Weingarten formulae give
	\begin{align}\label{position-correction-identities}
	\nabla_kv_i&=\delta_{ki}+\langle w,A_{ki}\rangle,\\
	\nabla_k^\perp w&=-\sum_pv_pA_{kp},\nonumber\\
	(\mathcal Lv)_i&=\sum_{k,p}\langle A_{ki},A_{kp}\rangle v_p.\nonumber
	\end{align}
For the last identity, both $\nabla_{\partial_t}v_i$ and $\Delta v_i$ contain $\langle w,\nabla_i^\perp H\rangle$, and these terms cancel. Since $|F|\le1$, it follows that
	\begin{align*}
	|v|+|w|&\le2,\qquad |\nabla v|\leq C(1+|A|),\qquad |\mathcal Lv|\leq C|A|^2.
	\end{align*}
Define the symmetric tensor of rank three with values in the normal bundle
	\begin{align}\label{B-explicit}
	B_{ijk}={}&v_i\left(\pi_{jk}A_{vv} -\sum_{p,q}\pi_j{}^p\pi_k{}^qA_{pq}\right)+v_j\left(\pi_{ik}A_{vv} -\sum_{p,q}\pi_i{}^p\pi_k{}^qA_{pq}\right)\\
	&+v_k\left(\pi_{ij}A_{vv} -\sum_{p,q}\pi_i{}^p\pi_j{}^qA_{pq}\right) +(2H-nA_{vv})v_iv_jv_k,\nonumber
	\end{align}
and set $Z=\nabla A-B$. Euclidean Codazzi implies that $Z$ is symmetric in all tangent indices.
\begin{lemma}\label{lem-B-construction}
The tensor defined by \eqref{B-explicit} satisfies
	\begin{align}\label{B-interior-bounds}
	\begin{split}
    |B|&\leq C|A|,\\
	|\nabla B|&\leq C|\nabla A|+C(1+|A|)|A|,\\
	|\mathcal LB|&\leq C(1+|A|)|\nabla A|+C(1+|A|)^3.
    \end{split}
	\end{align}
At the boundary its components are
	\begin{align}\label{Bequations}
    \begin{split}
	B_{\mu ab}&=\delta_{ab}A_{\mu\mu}-A_{ab},\qquad B_{\mu\mu\mu}=2H-nA_{\mu\mu},\qquad B_{abc}=B_{a\mu\mu}=0.
    \end{split}
	\end{align}
Moreover, at a boundary point in a frame parallel in the conormal direction,
	\begin{align}\label{B-conormal-derivatives}
	(\nabla_\mu B)_{abc}=(\nabla_\mu B)_{a\mu\mu}=0, \qquad |\nabla_\mu B|\leq C|A|.
	\end{align}
\end{lemma}
\begin{proof}
Write $B=\mathcal P(v)A$, where \eqref{B-explicit} defines the polynomial coefficient $\mathcal P$. The coefficients and their derivatives with respect to $v$ are bounded for $|v|\le1$. The first two estimates in \eqref{B-interior-bounds} follow from \eqref{position-correction-identities}. The product rule gives
	\begin{align*}
	\mathcal LB&=\mathcal P(v)\mathcal LA +D\mathcal P(v)[\mathcal Lv]A -\sum_kD^2\mathcal P(v)[\nabla_kv,\nabla_kv]A -2\sum_kD\mathcal P(v)[\nabla_kv]\nabla_kA.
	\end{align*}
The Euclidean evolution equation for $A$ gives $|\mathcal LA|\leq C|A|^3$, see \cite[Section 2]{LN}. This proves the last estimate in \eqref{B-interior-bounds}.

At the boundary, $v=\mu$, $w=0$, and $A_{a\mu}=0$. Thus \eqref{Bequations} follows by substitution. Also $\nabla_\mu v=\mu$ and $\nabla_\mu\pi=-2\mu\otimes\mu$ there. Differentiating \eqref{B-explicit}, every contribution to the components with indices $abc$ or $a\mu\mu$ vanishes: the possible remaining factors are $v_a$, $\nabla_\mu v_a$, or $A_{a\mu}$. For the other components only $(\nabla_\mu A)_{ab}$, $(\nabla_\mu A)_{\mu\mu}$ and $\nabla_\mu^\perp H$ occur. Their bounds follow from \eqref{eqn_nablamuAab}, \eqref{eqn_nablamuAmumu} and Lemma \ref{nablaH}. This proves \eqref{B-conormal-derivatives}.
\end{proof}

\begin{lemma}\label{lem-Z-boundary}
The tensor $Z=\nabla A-B$ satisfies the boundary conditions
	\begin{align}
	Z_{\mu ab}&=Z_{\mu\mu\mu}=0\label{Zids},\\
	(\nabla_\mu Z)_{abc}&=\delta_{bc}Z_{a\mu\mu}+\delta_{ab}Z_{c\mu\mu} +\delta_{ac}Z_{b\mu\mu}-2Z_{abc}\label{Zabc},\\
	(\nabla_\mu Z)_{a\mu\mu}&=2\sum_b Z_{abb}-(n+1)Z_{a\mu\mu}\label{Zmumumu}.
	\end{align}
Consequently,
	\begin{align}\label{partialmuZ}
	\partial_\mu|Z|^2&\leq C|Z|^2.
	\end{align}
In the interior,
	\begin{align}\label{Z-tensor-evolution}
	|\mathcal LZ|&\leq C(1+|A|)^2|Z|+C(1+|A|)^3,
	\end{align}
and hence
	\begin{align}\label{evolofZ}
	(\partial_t-\Delta)|Z|^2&\le-2|\nabla Z|^2+C(1+|A|)^2|Z|^2+C(1+|A|)^6.
	\end{align}
\end{lemma}
\begin{proof}
The identities \eqref{Zids} follow from \eqref{Bequations} and the spherical boundary identities. Fix a boundary frame with $\nabla^{\partial M}_{e_a}e_b=0$ and $\nabla_\mu e_i=0$ at the point. Then $\nabla_a\mu=e_a$ and $\nabla_a e_b=-\delta_{ab}\mu$. Differentiating \eqref{eqn_nablamuAab} tangentially gives
	\begin{align*}
	(\nabla^2A)_{a\mu bc}&=\delta_{bc}Z_{a\mu\mu}+\delta_{ab}Z_{c\mu\mu} +\delta_{ac}Z_{b\mu\mu}-2Z_{abc}.
	\end{align*}
Similarly, differentiation of \eqref{eqn_nablamuAmumu} gives
	\begin{align*}
	(\nabla^2A)_{a\mu\mu\mu}&=2\sum_bZ_{abb}-(n+1)Z_{a\mu\mu}.
	\end{align*}
The Ricci equation gives $R^\perp(\mu,e_a)=0$. The Gauss commutator terms acting on $A_{bc}$ or $A_{\mu\mu}$ contain a mixed component $A_{a\mu}$ and vanish. We may therefore interchange the derivative indices $a,\mu$ in these two formulae. Now use \eqref{B-conormal-derivatives} to obtain \eqref{Zabc} and \eqref{Zmumumu}. Only the components with zero or two conormal indices contribute to $\partial_\mu|Z|^2$, so \eqref{partialmuZ} follows.

Commuting the Euclidean evolution equation for $A$ gives $\mathcal L\nabla A=A*A*\nabla A$, where $*$ denotes contractions with the metrics. Combining this with Lemma \ref{lem-B-construction} and $\nabla A=Z+B$ proves \eqref{Z-tensor-evolution}. Taking the inner product with $2Z$ and using Young's inequality proves \eqref{evolofZ}. No derivative of $Z$ is present in the error term, so the full Bochner term is retained.
\end{proof}

\subsection{The first derivative estimate}

\begin{theorem}\label{thm_gradient}
Let $F:M^n\times[0,T)\to\overline B_R^{n+m}$, $n\ge4$, $m\ge2$, be a smooth compact MCF with free boundary on $\mathbb{S}_R^{n+m-1}$ and outward conormal $\mu=F/R$. Suppose initially
	\begin{align}\label{initialpinchinggradient}
	|A|^2&<\frac4{3n}|H|^2.
	\end{align}
Then, for a constant $C=C(n,m,M_0,R)$,
	\begin{align}\label{gradientestimate}
	|\nabla A|^2&\leq C\left(|A|^2+R^{-2}\right)^2, \quad\text{ on }M\times[0,T).
	\end{align}
\end{theorem}
\begin{proof}
Rescale to $R=1$. Lemma \ref{lem_fb-additive-gap} gives $\frac{5}{4n}<c_0<\frac{4}{3n}$ and $a>0$ such that \eqref{eq_fb-additive-gap} holds. In particular $|H|>0$ throughout the smooth flow, so the principal normal decomposition used in the reaction estimate is defined everywhere. Choose $c_0<c_1<\frac{4}{3n}$ and set $g=c_1|H|^2-|A|^2+1$. The preserved gap gives
	\begin{align}\label{gineq}
	1&\leq g,\qquad C^{-1}(|A|^2+1)\leq g\leq C(|A|^2+1).
	\end{align}
The Andrews--Baker reaction estimate \cite{AndrewsBaker} and the Kato inequality give
	\begin{align*}
	\frac12(\partial_t-\Delta)(g-1)&\geq (|h|^2+3|A^-|^2)(g-1) +\frac{5n-8}{9n}|\nabla A|^2.
	\end{align*}
In particular, for a constant $\eta>0$,
	\begin{align}\label{evolg}
	(\partial_t-\Delta)g&\ge\eta|\nabla A|^2\ge0.
	\end{align}
Use the explicit $B,Z$ constructed above, and write $u=|Z|^2$. Lemmas \ref{lem-B-construction} and \ref{lem-Z-boundary} imply
	\begin{align*}
	(\partial_t-\Delta)u&\le-2|\nabla Z|^2+Cgu+Cg^3,\\
	(\partial_t-\Delta)g&\ge\eta' u-Cg,
	\end{align*}
for some $\eta'>0$, since $|\nabla A|^2\ge\tfrac12|Z|^2-C|A|^2$. For $f=u/g^2$, the product rule and $|\nabla u|\le2\sqrt{u}|\nabla Z|$ give
	\begin{align}\label{ineqforf}
	(\partial_t-\Delta)f&\le\frac2g\langle\nabla f,\nabla g\rangle +Cg(f+1)-c g f^2
	\end{align}
with $c>0$. To see that no uncontrolled gradient term was discarded, the terms left after keeping the drift $2\nabla g/g$ are bounded by
	\begin{align*}
	-\frac2{g^2}|\nabla Z|^2 +\frac4{g^3}\sqrt u\,|\nabla Z||\nabla g| -\frac{2u}{g^4}|\nabla g|^2&=-\frac2{g^2} \left(|\nabla Z|-\frac{\sqrt u}{g}|\nabla g|\right)^2\le0.
	\end{align*}
The boundary pinching identity \eqref{evolQ} gives
	\begin{align}\label{partialmug}
	\partial_\mu g&\ge-2(g-1)\ge-2g.
	\end{align}
Together with the homogeneous bound \eqref{partialmuZ}, this yields
	\begin{align}\label{partialmuf2}
	\partial_\mu f&\leq Cf.
	\end{align}
Set $\rho=(1-|F|^2)/2$. Thus $0\le\rho\le\frac{1}{2}$, $|\nabla\rho|\le1$, $(\partial_t-\Delta)\rho=n$, and $\partial_\mu\rho=-1$. Choose a fixed $\lambda>C$ in \eqref{partialmuf2} and put $\widetilde f=e^{\lambda\rho}(f+1)$. Its outward derivative is strictly negative, so a spatial maximum occurs in the interior. At such a maximum, $\nabla f=-\lambda(f+1)\nabla\rho$. Before using the gradient estimate being proved, we have
	\begin{align*}
	|\nabla g|&\leq C|A||\nabla A| \leq C\bigl(g^{3/2}\sqrt f+g\bigr),
	\end{align*}
by $\nabla A=Z+B$ and \eqref{B-interior-bounds}. Consequently the drift and weight terms in \eqref{ineqforf} are absorbed by $\tfrac c2g f^2+C_\lambda g(f+1)$ and
	\begin{align*}
	(\partial_t-\Delta)\widetilde f&\leq e^{\lambda\rho} \left[C_\lambda g(f+1)-\frac c2g f^2\right].
	\end{align*}
The maximum principle therefore bounds $f$ in terms of its initial maximum and the fixed constants. We obtain
	\begin{align}\label{Zestimate}
	|Z|^2&\leq Cg^2, \qquad |\nabla A|^2\le2|Z|^2+2|B|^2\leq Cg^2.
	\end{align}
Rescaling proves \eqref{gradientestimate}.
\end{proof}

\subsection{The second derivative correction}

Continue with $R=1$, fix $c_0,c_1,g$ as in the preceding proof and use Theorem \ref{thm_gradient}. In particular,
	\begin{align}\label{first-grad-use}
	|\nabla A|+|Z|&\leq Cg,\qquad |\nabla g|\leq Cg^{3/2},\qquad |\nabla B|\leq Cg.
	\end{align}
We now prescribe every component of $\widetilde B$ explicitly. Define the tensors with values in the normal bundle
	\begin{align*}
	Z^\top_{ijk}&=\sum_{p,q,r}\pi_i{}^p\pi_j{}^q\pi_k{}^rZ_{pqr},\qquad \zeta_i=\sum_{p,q,r}\pi_i{}^p v_qv_rZ_{pqr},\qquad \tau_i=\sum_{p,q,r}\pi_i{}^p\pi_{qr}Z_{pqr}.
	\end{align*}
Set
	\begin{align}\label{Btilde-explicit}
	\hspace*{1cm}\widetilde B_{ijkl}&=v_i\Bigl[ \pi_{jk}\zeta_l+\pi_{jl}\zeta_k+\pi_{kl}\zeta_j-2Z^\top_{jkl}+v_jv_k\bigl(2\tau_l-(n+1)\zeta_l\bigr) \\
	&+v_jv_l\bigl(2\tau_k-(n+1)\zeta_k\bigr)+v_kv_l\bigl(2\tau_j-(n+1)\zeta_j\bigr)\Bigr]+v_j\bigl[-Z^\top_{ikl}+\pi_{ik}\zeta_l+\pi_{il}\zeta_k\bigr]\nonumber\\
	&+v_k\bigl[-Z^\top_{ijl}+\pi_{ij}\zeta_l+\pi_{il}\zeta_j\bigr]+v_l\bigl[-Z^\top_{ijk}+\pi_{ij}\zeta_k+\pi_{ik}\zeta_j\bigr] -3v_jv_kv_l\zeta_i,\nonumber
	\end{align}
and define $Y_{ijkl}=(\nabla_i Z)_{jkl}-\widetilde B_{ijkl}$. The tensor $\widetilde B$ is symmetric in its last three indices and is a bounded polynomial coefficient in $v$ applied linearly to $Z$. Thus it is defined on all of $M_t$, without a collar construction.

\begin{lemma}\label{lem-Y-corrections}
The correction \eqref{Btilde-explicit} satisfies
	\begin{align}\label{eq6}
	|\widetilde B|&\leq C|Z|\leq Cg, \qquad |\nabla^2A|^2\leq C(|Y|^2+g^2).
	\end{align}
At the free boundary,
	\begin{align}\label{eq7}
	\partial_\mu|Y|^2&\leq C(|Y|^2+g^3).
	\end{align}
In the interior,
	\begin{align}\label{eq8}
	(\partial_t-\Delta)|Y|^2&\le-2|\nabla Y|^2+Cg|Y|^2+Cg^4.
	\end{align}
All constants depend on the preceding first derivative estimate and not on any bound for $\nabla^2A$ or $\nabla^3A$.
\end{lemma}
\begin{proof}
The first bound in \eqref{eq6} follows directly from \eqref{Btilde-explicit}. Since $\nabla^2A=\nabla Z+\nabla B=Y+\widetilde B+\nabla B$, \eqref{first-grad-use} proves the second bound.

\smallskip\noindent\emph{Boundary components.} At the boundary, \eqref{Btilde-explicit} has the following nonzero component prescriptions, together with symmetry in the last three indices:
	\begin{align}\label{Btilde-boundary-components}
	\widetilde B_{\mu abc}&=\delta_{bc}Z_{a\mu\mu}+\delta_{ab}Z_{c\mu\mu} +\delta_{ac}Z_{b\mu\mu}-2Z_{abc},\\
	\widetilde B_{\mu a\mu\mu}&=2\sum_bZ_{abb}-(n+1)Z_{a\mu\mu},\nonumber\\
	\widetilde B_{a\mu bc}&=-Z_{abc}+\delta_{ab}Z_{c\mu\mu}+\delta_{ac}Z_{b\mu\mu},\nonumber\\
	\widetilde B_{a\mu\mu\mu}&=-3Z_{a\mu\mu}.\nonumber
	\end{align}
Every component with an even number of conormal indices is zero. The first two lines agree with \eqref{Zabc} and \eqref{Zmumumu}. The last two follow by differentiating \eqref{Zids} tangentially:
	\begin{align*}
	(\nabla_aZ)_{\mu bc}&=-Z_{abc}+\delta_{ab}Z_{c\mu\mu}+\delta_{ac}Z_{b\mu\mu}, \qquad (\nabla_aZ)_{\mu\mu\mu}=-3Z_{a\mu\mu}.
	\end{align*}
It follows that every component of $Y$ with an odd number of conormal indices vanishes at the boundary.

Let $a,b,c,d$ denote tangential indices along $\partial M_t$. Work in an adapted orthonormal frame satisfying
	\begin{align*}
	\nabla_a\mu&=e_a,\qquad \nabla_a e_b=-\delta_{ab}\mu,\qquad \nabla_\mu e_i=0
	\end{align*}
at the boundary point. Recall from \eqref{Zids} that $Z_{ab\mu}=Z_{\mu\mu\mu}=0$. Since every component of $Y$ with one or three conormal indices vanishes, only
	\begin{align*}
	Y_{abcd},\qquad Y_{ab\mu\mu},\qquad Y_{\mu ab\mu},\qquad Y_{\mu\mu\mu\mu}
	\end{align*}
can contribute to $\partial_\mu|Y|^2$, up to symmetry in the last three indices.

The conormal derivatives of the corresponding components of $\widetilde B$ vanish. Indeed, differentiating \eqref{Btilde-explicit} and using $v=\mu$ and $\nabla_\mu v=\mu$ gives
	\begin{align}\label{Btilde-conormal-even}
	(\nabla_\mu\widetilde B)_{abcd}&=0,\\
	(\nabla_\mu\widetilde B)_{ab\mu\mu}&=4Z_{ab\mu}-4\delta_{ab}Z_{\mu\mu\mu}=0,\nonumber\\
	(\nabla_\mu\widetilde B)_{\mu ab\mu}&=6Z_{ab\mu}-2\delta_{ab}Z_{\mu\mu\mu}=0,\nonumber\\
	(\nabla_\mu\widetilde B)_{\mu\mu\mu\mu}&=-12\sum_cZ_{cc\mu}+6(n+2)Z_{\mu\mu\mu}=0.\nonumber
	\end{align}
Terms containing $\nabla_\mu Z$ have zero coefficients in these components. The remaining terms vanish by \eqref{Zids}. Consequently, for the four component types above, differentiating $Y=\nabla Z-\widetilde B$ amounts to estimating the corresponding second derivatives of $Z$.

For $Y_{abcd}$ and $Y_{ab\mu\mu}$, differentiate \eqref{Zabc} and \eqref{Zmumumu} tangentially and interchange the two derivatives. The normal curvature terms vanish because
	\begin{align*}
	R^\perp(\mu,e_a)&=0.
	\end{align*}
The tangent curvature terms also vanish: the Gauss equation and $A_{a\mu}=0$ imply that $R(\mu,e_a)$ interchanges tangential and conormal directions. Its action on $Z_{bcd}$ or $Z_{b\mu\mu}$ therefore produces a component with one or three conormal indices, which is zero. Differentiating the adapted frame introduces only bounded multiples of $Z$ and $\nabla Z$. Together with \eqref{Btilde-conormal-even}, this gives
	\begin{align*}
	|(\nabla_\mu Y)_{abcd}|+|(\nabla_\mu Y)_{ab\mu\mu}|&\leq C\bigl(|\nabla Z|+|Z|\bigr).
	\end{align*}

The remaining components require two conormal derivatives of $Z$, which we obtain from its evolution equation. First, at the boundary,
	\begin{align}\label{Z-odd-evolution}
	|(\mathcal LZ)_{ab\mu}|+|(\mathcal LZ)_{\mu\mu\mu}|&\leq C\bigl(|A|^3+|A|\bigr).
	\end{align}
To prove this, expand
	\begin{align*}
	\mathcal LZ&=\mathcal L\nabla A-\mathcal LB
	\end{align*}
using $\mathcal L\nabla A=A*A*\nabla A$ and the explicit formula \eqref{B-explicit}. After using $A_{a\mu}=0$, every term involving $\nabla A$ in the $ab\mu$ and $\mu\mu\mu$ components contains only derivatives of the following types:
	\begin{align*}
	(\nabla_aA)_{b\mu}&=\delta_{ab}A_{\mu\mu}-A_{ab},\\
	(\nabla_\mu A)_{ab}&=\delta_{ab}A_{\mu\mu}-A_{ab},\\
	(\nabla_\mu A)_{\mu\mu}&=2H-nA_{\mu\mu}.
	\end{align*}
These are bounded by $C|A|$. The remaining terms are cubic or linear in $A$, by \eqref{position-correction-identities} and $\nabla_kv=e_k$ at the boundary. This proves \eqref{Z-odd-evolution}.

To recover the two conormal derivatives, differentiate the tangential identities obtained from $Z_{ab\mu}=Z_{\mu\mu\mu}=0$ once more and sum over $c$. Including the derivatives of the adapted frame gives
	\begin{align*}
	\sum_c(\nabla^2Z)_{ccab\mu} ={}&2(\nabla_aZ)_{b\mu\mu}+2(\nabla_bZ)_{a\mu\mu}-2\sum_c(\nabla_cZ)_{abc} +(n-1)(\nabla_\mu Z)_{ab\mu},\\
	\sum_c(\nabla^2Z)_{cc\mu\mu\mu} ={}&-6\sum_c(\nabla_cZ)_{c\mu\mu} +(n-1)(\nabla_\mu Z)_{\mu\mu\mu}.
	\end{align*}
The terms with coefficient $n-1$ arise from the Hessian correction $-\nabla_{\nabla_c e_c}Z=\nabla_\mu Z$ for each tangential direction $c$.

The time connection preserves the tangential and conormal subspaces, since $\mu=F$ and $\partial_tF=H$ is normal. The vanishing boundary components therefore satisfy
	\begin{align*}
	(\nabla_{\partial_t}Z)_{ab\mu}&=(\nabla_{\partial_t}Z)_{\mu\mu\mu}=0.
	\end{align*}
Using $\mathcal LZ=\nabla_{\partial_t}Z-\Delta Z$ and \eqref{Btilde-conormal-even}, we obtain the double conormal derivatives explicitly:
	\begin{align*}
	(\nabla_\mu Y)_{\mu ab\mu} ={}&-(\mathcal LZ)_{ab\mu} -2(\nabla_aZ)_{b\mu\mu}-2(\nabla_bZ)_{a\mu\mu}\\
	&+2\sum_c(\nabla_cZ)_{abc} -(n-1)(\nabla_\mu Z)_{ab\mu},\\
	(\nabla_\mu Y)_{\mu\mu\mu\mu} ={}&-(\mathcal LZ)_{\mu\mu\mu} +6\sum_c(\nabla_cZ)_{c\mu\mu}-(n-1)(\nabla_\mu Z)_{\mu\mu\mu}.
	\end{align*}
Every term on the right is now controlled by \eqref{Z-odd-evolution} or a first derivative of $Z$. Hence
	\begin{align*}
	|(\nabla_\mu Y)_{\mu ab\mu}|+|(\nabla_\mu Y)_{\mu\mu\mu\mu}|&\leq C\bigl(|\nabla Z|+|A|^3+|A|\bigr).
	\end{align*}

We have controlled every component that contributes to $\partial_\mu|Y|^2$. By \eqref{eq6}, \eqref{first-grad-use} and the comparison between $g$ and $1+|A|^2$,
	\begin{align*}
	\nabla Z&=Y+\widetilde B,\qquad |\widetilde B|+|Z|\leq Cg,\qquad 1+|A|^2\leq Cg.
	\end{align*}
Consequently,
	\begin{align*}
	\partial_\mu|Y|^2&=2\sum_{i,j,k,l} \bigl\langle(\nabla_\mu Y)_{ijkl},Y_{ijkl}\bigr\rangle\\
	&\leq C|Y|\bigl(|\nabla Z|+|Z|+|A|^3+|A|\bigr)\\
	&\leq C|Y|\bigl(|Y|+g^{3/2}\bigr) \leq C\bigl(|Y|^2+g^3\bigr).
	\end{align*}
This proves \eqref{eq7} without assuming a bound for $\nabla Y$.

\smallskip\noindent\emph{Interior evolution.} The first derivative estimate and \eqref{position-correction-identities} give $|\nabla^2v|\leq Cg$. Differentiating the explicit product formula for $\mathcal LB$, and using $\mathcal L\nabla A=A*A*\nabla A$, gives
	\begin{align*}
	|\nabla\mathcal LZ|&\leq Cg|\nabla Z|+Cg^{5/2}.
	\end{align*}
Here $\nabla^2A=\nabla Z+\nabla B$, while every term containing two first derivatives of $A$ is bounded by $C|A||\nabla A|^2\leq Cg^{5/2}$. The connection commutator satisfies
	\begin{align*}
	|[\mathcal L,\nabla]Z|&\leq C\bigl(|A|^2|\nabla Z|+|A||\nabla A||Z|\bigr).
	\end{align*}
Write \eqref{Btilde-explicit} as $\widetilde B=\mathcal C(v)Z$. Its product rule, \eqref{Z-tensor-evolution}, and \eqref{first-grad-use} imply
	\begin{align*}
	|\mathcal L\widetilde B|&\leq C\sqrt g\,|\nabla Z|+Cg^2.
	\end{align*}
Thus
	\begin{align*}
	|\mathcal LY|&\leq Cg|\nabla Z|+Cg^{5/2} \leq Cg|Y|+Cg^{5/2}.
	\end{align*}
Taking the inner product with $2Y$ and using $2Cg^{5/2}|Y|\leq Cg|Y|^2+Cg^4$ proves \eqref{eq8}. No bound for $\nabla Y$ was used in these error estimates.
\end{proof}

\subsection{The second derivative estimate}
\begin{theorem}\label{thm_second_gradient}
Let $F:M^n\times[0,T)\to\overline B_R^{n+m}$, $n\ge4$, $m\ge2$, be a smooth compact MCF with free boundary on $\mathbb{S}_R^{n+m-1}$, satisfying \eqref{initialpinchinggradient} at $t=0$. Then, for $C=C(n,m,M_0,R)$,
	\begin{align*}
	|\nabla^2A|^2&\leq C\left(|A|^2+R^{-2}\right)^3, \quad\text{ on }M\times[0,T).
	\end{align*}
\end{theorem}
\begin{proof}
Rescale to $R=1$ and use the same $c_0,c_1$ and $g$ as in the first derivative estimate. The preserved gap and the positive lower bound for $|H|$ follow from Lemma \ref{lem_fb-additive-gap}. Set $u=|Z|^2$, $f=u/g^2$, and $q=|Y|^2/g^3$. By Theorem \ref{thm_gradient}, $f\leq C$. Since $\nabla Z=Y+\widetilde B$, \eqref{evolofZ}, \eqref{gineq} and \eqref{eq6} give
	\begin{align}\label{eq9}
	(\partial_t-\Delta)u&\le-|\nabla Z|^2-\tfrac12|Y|^2+Cg^3.
	\end{align}
We retain the exact quotient drifts. Using \eqref{evolg}, \eqref{first-grad-use} and \eqref{eq9},
	\begin{align}\label{eq11}
	(\partial_t-\Delta)f&\le\frac4g\langle\nabla f,\nabla g\rangle-cgq+Cg
	\end{align}
for $c>0$. Indeed, the remaining gradient term from the denominator is $2f|\nabla g|^2/g^2\leq Cg$. Similarly \eqref{eq8} gives
	\begin{align}\label{eq12}
	(\partial_t-\Delta)q&\le\frac6g\langle\nabla q,\nabla g\rangle+Cg(q+1),
	\end{align}
since the corresponding term is $6q|\nabla g|^2/g^2\leq Cgq$. The different coefficients in \eqref{eq11} and \eqref{eq12} must be accounted for. From \eqref{first-grad-use} and \eqref{eq6},
	\begin{align*}
	|\nabla f|&\leq C\sqrt g(\sqrt q+1),\qquad \frac{|\nabla f||\nabla g|}{g}\leq Cg(\sqrt q+1).
	\end{align*}
For $N\ge1$, set $\Phi=q+Nf+1$. Adding \eqref{eq12} to $N$ times \eqref{eq11}, the drift is exactly
	\begin{align*}
	\frac6g\langle\nabla\Phi,\nabla g\rangle -\frac{2N}{g}\langle\nabla f,\nabla g\rangle.
	\end{align*}
Young's inequality therefore yields
	\begin{align}\label{Phi-interior-corrected}
	(\partial_t-\Delta)\Phi&\le\frac6g\langle\nabla\Phi,\nabla g\rangle -\left(\frac{cN}{2}-C\right)gq+C_Ng.
	\end{align}

At the boundary, \eqref{eq7} and \eqref{partialmug} imply $\partial_\mu q\leq C(q+1)$, whereas \eqref{partialmuf2} is homogeneous. Consequently
	\begin{align}\label{eq15}
	\partial_\mu\Phi&\leq C\Phi,
	\end{align}
with $C$ independent of $N$. First choose a fixed $\lambda>C$ in \eqref{eq15}, and then set $\widetilde\Phi=e^{\lambda\rho}\Phi$, where $\rho=(1-|F|^2)/2$. Its outward derivative is strictly negative. At an interior spatial maximum, $\nabla\Phi=-\lambda\Phi\nabla\rho$. Using $|\nabla g|/g\leq C\sqrt g$, $g\ge1$ and $f\leq C$, all drift and weight terms in \eqref{Phi-interior-corrected} are bounded above by
	\begin{align*}
	C_\lambda g\Phi&\leq C_\lambda gq+C_{\lambda,N}g,
	\end{align*}
where $C_\lambda$ is independent of $N$. We now choose $N$ so large that $cN/2-C-C_\lambda\ge1$. At the maximum,
	\begin{align*}
	(\partial_t-\Delta)\widetilde\Phi&\leq e^{\lambda\rho}(-gq+C_{\lambda,N}g).
	\end{align*}
Since $q\geq e^{-\lambda/2}\widetilde\Phi-N\sup f-1$, this inequality bounds $\widetilde\Phi$, and hence $q$, in terms of its initial maximum and the fixed constants. Thus $|Y|^2\leq Cg^3$. Finally \eqref{eq6} and \eqref{gineq} imply $|\nabla^2A|^2\leq Cg^3\leq C(|A|^2+1)^3$. Rescaling proves the theorem. The quotient and boundary estimates above use only the positivity in \eqref{evolg} and the correction identities and apply for $n\geq4$.
\end{proof}

\section{Boundary Proof}\label{sec:boundary-proof}

\begin{corollary}[Derivative estimates under initial quadratic convexity]
Let $n\geq4$, $m\geq2$ and $\kappa=14/5$. Suppose a smooth compact MCF with free boundary on $\mathbb{S}_R^{n+m-1}$ satisfies $F_\e>0$ initially, for some fixed $\e>0$. If $4\leq n\leq7$, assume in addition that \eqref{initialpinchinggradient} holds initially. Then, for $C=C(n,m,M_0,R)$,
	\begin{align*}
	|\nabla A|^2&\leq C\left(|A|^2+R^{-2}\right)^2, \qquad |\nabla^2A|^2\leq C\left(|A|^2+R^{-2}\right)^3
	\end{align*}
throughout the smooth flow.
\end{corollary}
\begin{proof}
Lemma \ref{lem_modified-quadratic-pinching}, specifically \eqref{eq_modified-strict-coarse} at $t=0$, gives \eqref{initialpinchinggradient}, for $n\ge8$. For $4\leq n\le7$, it holds by assumption. The conclusion follows directly from Theorems \ref{thm_gradient} and \ref{thm_second_gradient}.
\end{proof}

\subsection{Boundary derivative lemma}
The boundary step specific to quadratic convexity is a strictly positive outward conormal derivative at a possible first zero. We prove this for a minimum eigenvector tangent to the boundary and for the conormal eigenvector. Combined with Lemma \ref{lem_interior}, the estimate completes the preservation argument for spherical free boundary flow in arbitrary codimension.
\begin{lemma}[Boundary inequality]\label{lem_boundary}
Let $F:M^n\times[0,T)\to\overline B_R^{n+m}$, $n\ge4$, be a smooth MCF with free boundary on $\mathbb S_R^{n+m-1}$ and outward conormal $\mu=F/R$. Set
	\begin{align*}
	c_\partial(\kappa)&=\frac{-4\kappa^2+12\kappa-1}{8\kappa}.
	\end{align*}
Suppose $\kappa>0$ and $c_\partial(\kappa)>0$. There exists $\e_0=\e_0(n,\kappa)>0$ such that, for every $0<\e\le\e_0$, the following holds. Suppose that $F_\e=0$ at a boundary point with $|H|>0$. At that point, choose a minimum eigenvector $e$ tangent to $\partial M$ or equal to $\mu$, and extend it as a unit vector field with $\nabla_\mu e=0$ at that point. Then
	\begin{align*}
	\partial_\mu F_{\e,e}&\ge\frac{c_\partial(\kappa)}{2nR}|H|^2>0.
	\end{align*}
\end{lemma}
\begin{proof}
By Lemma \ref{lem_modified-principal-convexity}, we have \eqref{rangeoflambda1}, and $F_{\e,e}$ is a smooth upper support at the chosen minimum eigenvector. The equation at the zero can also be written as
	\begin{align}\label{eqn_boundary-zero}
	\kappa|A|^2&=|H|\lambda_1 -\frac{n-1}{2}\bigl(\lambda_1^2+|A^-|^2\bigr) +\kappa\left(\frac1{n-1}-\e\right)|H|^2.
	\end{align}
By Lemma \ref{nablaH}, together with \eqref{eq_quadratic-boundary}, \eqref{partialmuA} and \eqref{norotation}, we can compute the conormal derivative of $F_{\e,e}$. We distinguish the cases in which a minimum eigenvector is tangential or conormal.
Since $h_{a\mu}=0$, a minimum eigenvector can be chosen in the boundary tangent space or the conormal line. For a tangential minimum eigenvector, $R\nabla_\mu h(e,e)=h_{\mu\mu}-\lambda_1$ and similarly in the conormal case, $R(\nabla_\mu h)_{\mu\mu}=2|H|-n\lambda_1$. Differentiating \eqref{eqn_F} step by step, we have
	\begin{align*}
	R \partial_\mu(|H| \lambda_1)&=|H| \lambda_1+|H|\left(h_{\mu \mu}-\lambda_1\right)=|H| h_{\mu \mu}.\\
	R \partial_\mu\left[-\frac{n-1}{2} \lambda_1^2\right]&=-(n-1) \lambda_1\left(h_{\mu \mu}-\lambda_1\right)\\
	R \partial_\mu\left[-\frac{n-1}{2}\left|A^-\right|^2\right]&=(n-1)\left|A^-\right|^2+n(n-1)\left|A_{\mu \mu}^-\right|^2\\
		\hspace*{-1.7cm}R \partial_\mu\left\{\kappa\left[\left(\frac{1}{n-1}-\e\right) |H|^2-|A|^2\right]\right\}&= \kappa\Big[2\left(\frac{1}{n-1}-\e\right) |H|^2 +2|A|^2-6 |H|h_{\mu \mu}\\
	&+2 n h_{\mu \mu}^2+2 n\left|A_{\mu \mu}^-\right|^2\Big].
	\end{align*}
Eliminating $\kappa|A|^2$ by \eqref{eqn_boundary-zero} gives
	\begin{align*}
	R\partial_\mu F_{\e,e}&=2\kappa n h_{\mu\mu}^2-\bigl((6\kappa-1)|H|+(n-1)\lambda_1\bigr)h_{\mu\mu}+4\kappa\left(\frac{1}{n-1}-\e\right)|H|^2\\
	&+2|H|\lambda_1+n(n-1+2\kappa)|A^-_{\mu\mu}|^2.
	\end{align*}
For fixed $h_{\mu\mu}$ this is affine in $\lambda_1\in[0,|H|/n]$. At the two endpoints, square completion gives the respective coefficients
	\begin{align*}
	\frac{4\kappa}{n-1}-\frac{(6\kappa-1)^2}{8\kappa n}, \qquad \frac2n+\frac{4\kappa}{n-1}-\frac{(6\kappa-1/n)^2}{8\kappa n}
	\end{align*}
in front of $|H|^2$, before subtracting $4\kappa\e$. The second exceeds the first by
	\begin{align*}
	\frac{n+3}{2n^2}+\frac{n^2-1}{8\kappa n^3}&>0.
	\end{align*}
The first is at least $c_\partial(\kappa)/n$, where
	\begin{align*}
	c_\partial(\kappa)&=\frac{-4\kappa^2+12\kappa-1}{8\kappa}.
	\end{align*}
For a conormal minimum eigenvector, $h_{\mu\mu}=\lambda_1$ and $R\nabla_\mu h_{\mu\mu}=2|H|-n\lambda_1$.
This gives
	\begin{align*}
	R\partial_\mu F_{\e,e}&= \bigl((n-1)^2+2\kappa n\bigr)\lambda_1^2-(3n+6\kappa-5)|H|\lambda_1\\
	&+\left(2+4\kappa\left(\frac{1}{n-1}-\e\right)\right)|H|^2+n(n-1+2\kappa)|A^-_{\mu\mu}|^2.
	\end{align*}
The quadratic term is decreasing on $[0,|H|/n]$, since its derivative at the right endpoint, divided by $|H|$, is $-n+1+2/n-2\kappa<0$. Thus its lower bound is
	\begin{align*}
	\left(\frac3n+\frac1{n^2}+\frac{4\kappa}{n(n-1)}-4\kappa\e\right)|H|^2.
	\end{align*}
Since $c_\partial(\kappa)\le1$ for every $\kappa>0$, both cases give
	\begin{align*}
	R\partial_\mu F_{\e,e}&\ge\left(\frac{c_\partial(\kappa)}n-4\kappa\e\right)|H|^2.
	\end{align*}
The coefficient $c_\partial(\kappa)$ is positive precisely when
	\begin{align*}
	\frac32-\sqrt2&<\kappa<\frac32+\sqrt2.
	\end{align*}
Fix such a $\kappa$. Since $c_\partial(\kappa)/n>0$, we may choose $\e_0=\e_0(n,\kappa)>0$ sufficiently small that
	\begin{align*}
	\frac{c_\partial(\kappa)}n-4\kappa\e&\ge\frac{c_\partial(\kappa)}{2n}>0, \qquad\text{ for every } 0<\e\le\e_0.
	\end{align*}
The outward derivative is therefore at least $c_\partial(\kappa)|H|^2/(2nR)$, as required.
\end{proof}

\subsection{Preservation of quadratic convexity}
\begin{theorem}[Preservation of quadratic convexity for spherical free boundary flow]\label{theorem8.3}
\label{thm_free-boundary-preservation}
Let $F:M^n\times[0,T)\to\overline B_R^{n+m}$ be a smooth compact MCF with free boundary on $\mathbb{S}_R^{n+m-1}$, where $n\ge4$. Fix $\kappa=14/5$ and let $\e>0$ be sufficiently small. Suppose that $F_\e$, as defined in \eqref{eqn_F}, is strictly positive initially. If $4\leq n\leq7$, assume in addition that \eqref{initialpinchinggradient} holds initially. Then $F_\e>0$ throughout the smooth flow. Moreover, there exist constants $a>0$ and $\frac{1}{n}<c<\frac{4}{3n}$, such that
	\begin{align*}
	|A|^2+a&\leq c|H|^2,\qquad |H|^2\ge\frac{a}{c-1/n}>0
	\end{align*}
throughout the flow. In particular,
	\begin{align*}
	h&\ge\frac{\kappa\e}{1+2\kappa/(n-1)}|H|g.
	\end{align*}
\end{theorem}
\begin{proof}
For $n\ge8$, the strict initial quadratic pinching follows from \eqref{eq_modified-strict-coarse}, while for $4\leq n\le7$ it holds by assumption. Lemma \ref{lem_fb-additive-gap} therefore gives constants $a>0$ and $c<\frac{4}{3n}$ such that
	\begin{align*}
	|A|^2+a&\leq c|H|^2
	\end{align*}
throughout the smooth flow. With $c=c_0$ in that lemma, \eqref{eq_fb-additive-gap} gives the preserved additive gap and the positive lower bound for $|H|$. In particular the principal normal remains defined.

\smallskip\noindent\emph{Preservation of $F_\e>0$.} Since $c_\partial(14/5)=31/560>0$, take $\e>0$ sufficiently small for Lemma \ref{lem_boundary}. Compactness gives $\min_{M_0}F_\e/|H|^2>0$. Choose $\e'>\e$ sufficiently close that $F_{\e'}>0$ on $M_0$ and $\e'$ also satisfies its smallness restriction.

Fix $S<T$ and take $\sigma>0$ small enough that $\delta(t)=\e'-\sigma t>0$ on $[0,S]$. At a first zero of $F_{\delta(t)}$, Lemma \ref{lem_modified-principal-convexity} provides a smooth upper support $F_{\delta(t),e}$ with a local minimum equal to zero at the point of contact. If the contact is interior, Lemma \ref{lem_interior} applies with the coefficient held constant at its value at that time, since \eqref{eq_fb-additive-gap} implies the quadratic pinching. The time dependence of $\delta$ contributes $\kappa\sigma|H|^2$, so
	\begin{align*}
	(\partial_t-\Delta)F_{\delta(t),e}&\ge\kappa\sigma|H|^2>0,
	\end{align*}
contradicting the differential inequality at the first minimum. At a boundary contact, choose $e$ and its extension as in Lemma \ref{lem_boundary}. That lemma gives $\partial_\mu F_{\delta(t),e}>0$, whereas the outward derivative at the minimum is nonpositive. Letting $\sigma\to 0$ and taking $S<T$ arbitrary gives $F_{\e'}\ge0$ throughout the smooth flow. Hence
	\begin{align*}
	F_\e&=F_{\e'}+\kappa(\e'-\e)|H|^2 \ge\kappa(\e'-\e)|H|^2>0.
	\end{align*}
The estimate for $h$ follows from \eqref{rangeoflambda1}.
\end{proof}

\section{Free Boundary Convergence}\label{sec:fb-convergence}

The derivative estimates established above provide the compactness needed to analyse singularities. We use the compactness and continuation result in Appendix \ref{app:fb-compactness}, which adapts the graphical compactness theory for immersed submanifolds to spherical free boundary MCF. It gives smooth pointed convergence up to the boundary and continuation whenever the curvature and its first two covariant derivatives remain bounded. When the rescaled support spheres converge to a hyperplane, reflection extends the limiting flow smoothly across its boundary.

We combine known rigidity results for ancient flows with boundary compactness to prove contraction to a round hemisphere. Quantitative planarity then gives an integrable bound for the rotation of its ambient plane near extinction. This proves uniqueness of the limiting plane and smooth convergence of the whole rescaled flow up to the boundary.

\subsection{Classification of ancient limits}

\begin{lemma}\label{lem_fb-ancient-rigidity}
Let $n\ge4$ and let a nonflat complete connected ancient MCF have bounded curvature on compact time intervals. Suppose
	\begin{align*}
	|A|^2&\leq c|H|^2,\qquad h\ge\eta|H|g, \qquad c<\frac4{3n},\quad\eta>0,
	\end{align*}
where the condition on $h$ is imposed wherever $H\ne0$. If $n\ge5$, assume in addition that $A^-=0$ wherever $H\ne0$. If the flow has no boundary, it consists of shrinking round spheres in a fixed affine $(n+1)$-plane. If instead it has nonempty orthogonal free boundary on a hyperplane and lies on one side of it, the conclusion is shrinking round hemispheres.
\end{lemma}

\begin{proof}
In the boundary case, reflect across the support using Lemma \ref{lem_fb-smooth-compactness}. For $c<c_1<\frac{4}{3n}$, the quadratic gap $c_1|H|^2-|A|^2$ is a nonnegative supersolution by the Andrews--Baker reaction estimate \cite{AndrewsBaker}, as in \eqref{evolg}. The strong maximum principle and uniqueness for a flow that is totally geodesic at one time therefore give $|H|>0$ throughout the ancient flow. For $n=4$, we have $c<1/3$, and hence
	\begin{align*}
	|A|^2-\frac13|H|^2&\leq -\left(\frac13-c\right)|H|^2.
	\end{align*}
Taking $0<\delta<1/3-c$ gives $|A|^2-\frac13|H|^2<-\delta|H|^2$. At every fixed time, the limiting submanifold is complete with bounded, nonvanishing mean curvature, so \cite[Corollary 3.3]{LNAncient2021} implies that it is compact. The ancient rigidity theorem \cite[Theorem 1.1]{LNAncient2021} then shows that the flow consists of shrinking round spheres. Uniqueness of the flow starting from a round sphere implies that these spheres have a fixed centre and a fixed affine $(n+1)$-plane. In the boundary case, reflection symmetry gives shrinking round hemispheres. This proves the assertion for $n=4$. Henceforth, assume $n\ge5$. Since
	\begin{align*}
	|H|g-h&\ge(1-\sqrt c)|H|g,
	\end{align*}
Naff's result on codimension reduction \cite[Proposition 2.5]{Naff2022} applies. At each time, the nonflat limiting submanifold is a complete strictly convex hypersurface in an affine $(n+1)$-plane. Completeness and strict convexity imply that the immersion is embedded and bounds a convex region. Hamilton's compactness theorem \cite[Theorem 1.1]{HamiltonCompactness1994} then applies to $h\ge\eta|H|g$, so the domain is compact.

Uniqueness excludes a totally geodesic submanifold at an earlier time, and nonvanishing of the mean curvature persists at later times. The affine plane is fixed because its normal coordinate functions solve the heat equation. Extend forward by Huisken's convex flow \cite{Huisken1984} and shift its extinction time to zero. The pinching remains uniform, so Huisken-Sinestrari's ancient rigidity theorem \cite[Theorem 3.2]{HuiskenSinestrari2015} gives shrinking spheres. Reflection symmetry then gives hemispheres in the boundary case.
\end{proof}

\subsection{Convergence to a round hemisphere}

\begin{proof}[Proof of Theorem \ref{thm_fb-convergence}]
Theorem \ref{theorem8.3} gives
	\begin{align*}
	|A|^2+a&\leq c|H|^2,\qquad h\ge\eta|H|g, \qquad a,\eta>0,\quad c<\frac4{3n}.
	\end{align*}
In particular, $|H|$ has a positive lower bound.

\smallskip\noindent\emph{Singularity formation in finite time.} Write $V(t)=\operatorname{Area}(M_t)$. The Michael-Simon inequality \cite{MichaelSimon1973} with boundary term and the identity
	\begin{align*}
	R\operatorname{Area}(\partial M_t)&=nV(t)+\int_{M_t}\langle F,H\rangle
	\end{align*}
give
	\begin{align*}
	V(t)^{(n-1)/n}&\leq C\left(\int_{M_t}|H|+\operatorname{Area}(\partial M_t)\right) \leq C\int_{M_t}|H|.
	\end{align*}
Here we used $|F|\leq R$ and the positive lower bound for $|H|$. Since $V'=-\int_{M_t}|H|^2$, Cauchy-Schwarz yields
	\begin{align*}
	V'&\le-C^{-1}V^{(n-2)/n},
	\end{align*}
and hence
	\begin{align}\label{Testimate}
	T-t&\leq CV(t)^{2/n}.
	\end{align}
Thus $T<\infty$, and Lemma \ref{lem_fb-smooth-compactness}(iii) together with Theorems \ref{thm_gradient} and \ref{thm_second_gradient} implies unbounded curvature at $T$.

\smallskip\noindent\emph{Blowup and contraction to a point.} Use rescaling at times when the running maximum is attained, as in Lynch--Nguyen's blowup argument \cite[Section 5]{LN}, with Lemma \ref{lem_fb-smooth-compactness}. Specifically, choose
	\begin{align*}
	Q_j&=|H|(p_j,t_j) =\max_{M\times[0,t_j]}|H|\longrightarrow\infty
	\end{align*}
and set
	\begin{align*}
	F_j(p,\tau)&=Q_j\bigl(F(p,t_j+Q_j^{-2}\tau)-F(p_j,t_j)\bigr), \qquad -Q_j^2t_j\le\tau\le0.
	\end{align*}
The curvature and derivative estimates give a complete nonflat ancient limit on $(-\infty,0]$.
For $n\ge5$, the planarity estimate \eqref{fbestimate} gives $A^-_\infty=0$ and Lemma \ref{lem_fb-ancient-rigidity} classifies the ancient limit. For $n=4$, the inherited bound $c<1/3$ is sufficient to apply Lemma \ref{lem_fb-ancient-rigidity}, without assuming $A^-_\infty=0$. Thus, in either dimension range, that lemma classifies the ancient limit as a shrinking round sphere or hemisphere.

Because the limit is compact, its convergence maps are eventually defined on the whole limiting manifold. Their images are open and closed in the connected manifold $M$, so convergence is global. Since $\partial M\ne\varnothing$, the limit is a hemisphere. Its radius at time zero is $n$. Hence
	\begin{align*}
	Q_j^nV(t_j)&\longrightarrow\frac12|\mathbb S^n|n^n,
	\end{align*}
and \eqref{Testimate} gives
	\begin{align}\label{T-tj}
	(T-t_j)Q_j^2&\leq C.
	\end{align}
The latter bound implies Type I. Indeed, any sequence violating Type I can be replaced by earlier times when the running maximum is attained, for which $(T-t)\max|H|^2$ is at least as large. Consequently,
	\begin{align}\label{typeIbound}
	|A|^2&\le\frac C{T-t},\qquad |H|\le\frac C{\sqrt{T-t}}.
	\end{align}
The speed is integrable, and global convergence along the sequence of times when the running maximum is attained makes the unscaled diameters tend to zero. Thus the flow contracts to a single point $q\in\mathbb S_R^{n+m-1}$, with
	\begin{align}\label{eq_fb-extinction-position}
	|F(p,t)-q|&\leq C\sqrt{T-t}.
	\end{align}

\smallskip\noindent\emph{Limits at the extinction scale.} Given any $\widehat t_j\to T$, put $l_j=\sqrt{T-\widehat t_j}$ and
	\begin{align*}
	\widehat F_j(p,\tau)&=l_j^{-1}\bigl(F(p,T+l_j^2\tau)-q\bigr), \qquad -T/l_j^2\le\tau<0.
	\end{align*}
With a boundary basepoint at time $-1$, the Type I and derivative estimates give smooth subsequential convergence on every $[-S,-\delta]$, $S>1$, $0<\delta<1$. The limit has free boundary on $q^\perp$, lies in $\langle q,x\rangle\le0$, and by \eqref{eq_fb-extinction-position} satisfies
	\begin{align*}
	|F_\infty(p,\tau)|&\leq C\sqrt{-\tau}.
	\end{align*}
Its reflected double cannot be flat, since a complete totally geodesic immersion has unbounded image. When $n\geq5$, the planarity estimate again gives $A^-_\infty=0$. When $n=4$, the strict bound $c<1/3$ is inherited. Lemma \ref{lem_fb-ancient-rigidity} therefore gives shrinking hemispheres. The position bound fixes their centre at zero and extinction time at zero.

Every sequence consequently has a globally smoothly convergent subsequence whose limiting submanifold at time $-1$ is a hemisphere of radius $\sqrt{2n}$. In particular, for $r=\sqrt{T-t}$ and all sufficiently late times,
	\begin{align}\label{eq_fb-late-slice-bounds}
	r|H|&\geq b>0,\qquad \diam(M,g(t))\leq Cr.
	\end{align}

\smallskip\noindent\emph{Uniqueness of the limiting plane.} For $n=4$, every sequence of rescalings at the extinction scale has a globally smoothly convergent subsequence whose limiting submanifold at time $-1$ is a round hemisphere. The mean curvature of the limit is nonzero.
Thus
	\begin{align*}
	\lim_{t\to T}\sup_{M_t}\frac{|A|^2}{|H|^2}&=\frac14.
	\end{align*}
Choose
	\begin{align*}
	\frac9{32}&<c_0<c_1<\frac5{16}.
	\end{align*}
Because the limiting ratio is $1/4<c_0$, there is $t_0<T$ such that $|A|^2<c_0|H|^2$ at $t=t_0$. Proposition \ref{prop_fb-quadratic-preservation} applies because $c_0>\frac{9}{8n}=9/32$ and preserves this inequality on $[t_0,T)$. Hence,
	\begin{align*}
	Q&=c_1|H|^2-|A|^2\ge(c_1-c_0)|H|^2>0.
	\end{align*}
The local calculations in the second cases of \cite[Lemmas 4.6-4.9]{Naff2022} remain valid for $n=4$, where the coefficient in their pinching assumption is $3(n+1)/(2n(n+2))=5/16$. Since $c_1<5/16$, these estimates, together with \cite[Lemmas 4.4-4.5 and Section 4.3]{Naff2022}, give, for sufficiently small $\sigma>0$,
	\begin{align*}
	(\partial_t-\Delta)P_\sigma&\le\frac{2(1-\sigma)}{Q}\langle\nabla Q,\nabla P_\sigma\rangle,\qquad P_\sigma=\frac{|A^-|^2}{Q^{1-\sigma}}.
	\end{align*}
The boundary calculations \eqref{ineqQ} and \eqref{partialmuA} remain valid and give
	\begin{align*}
	\frac R2\partial_\mu P_\sigma&\le-\sigma P_\sigma-\frac{4|A^-_{\mu\mu}|^2}{Q^{1-\sigma}}.
	\end{align*}
The maximum principle therefore yields
	\begin{align*}
	|A^-|^2&\leq C|H|^{2-2\sigma},
	\end{align*}
on $M\times[t_0,T)$.

For $n\ge5$, the planarity estimate \eqref{fbestimate} holds throughout the flow, while for $n=4$ the preceding argument gives the same estimate on $[t_0,T)$. Thus, in either case, on the extinction scale we have $|\widetilde A^-|\leq Cr^\sigma$.

Lemma \ref{lem_fb-smooth-compactness}(i) supplies uniform higher derivative bounds on rescaled backward intervals, and \eqref{eq_fb-late-slice-bounds} allows differentiation of the decomposition along the principal normal. Interpolation between the $C^0$ decay and uniform $C^4$ bounds for the rescaled $A^-$, using bounded extension operators in charts on half balls at the boundary, gives, with $\alpha=\sigma/2>0$,
	\begin{align}\label{eq_fb-planarity-derivatives}
	|A^-|&\leq Cr^{-1+\sigma},\qquad |\nabla A^-|\leq Cr^{-2+\alpha},\qquad |\nabla^2A^-|\leq Cr^{-3+\alpha}.
	\end{align}
Set $w_i=\nabla_i^\perp\nu$, and let $\pi$ be normal projection onto $\nu^\perp$. Contracted Codazzi gives
	\begin{align}\label{eq_fb-codazzi-normal}
	\sum_i(|H|\delta_{ki}-h_{ki})w_i&=\sum_i\pi\bigl((\nabla_iA^-)_{ki}\bigr).
	\end{align}
Since $|H|g-h\ge(1-\sqrt c)|H|g$, \eqref{eq_fb-codazzi-normal} and its first derivative imply
	\begin{align*}
	|w|&\leq Cr^{-1+\alpha},\qquad |\nabla w|\leq Cr^{-2+\alpha}.
	\end{align*}
The normal evolution equation
	\begin{align*}
	\nabla_{\partial_t}^\perp\nu&=\pi(\Delta^\perp\nu) +2\sum_i(\nabla_i\log|H|)w_i +\sum_{i,j}h_{ij}A^-_{ij}
	\end{align*}
then gives $|\nabla_{\partial_t}^\perp\nu|\leq Cr^{-2+\alpha}$. Let $\Pi(p,t)$ denote ambient orthogonal projection onto
	\begin{align*}
	dF(T_pM)\oplus\mathbb R\nu.
	\end{align*}
Differentiating an orthonormal frame gives
	\begin{align}\label{eq_fb-rotation}
	|D\Pi|&\leq C\bigl(|A^-|+|w|\bigr) \leq Cr^{-1+\alpha},\\
	|\partial_t\Pi|&\leq C\bigl(|H||w|+|\nabla_{\partial_t}^\perp\nu|\bigr) \leq C(T-t)^{-1+\alpha/2}.\nonumber
	\end{align}
The time bound in \eqref{eq_fb-rotation} is integrable, while the spatial oscillation is at most
	\begin{align*}
	Cr^{-1+\alpha}\diam(M,g(t))&\leq Cr^\alpha.
	\end{align*}
Hence $\Pi$ converges uniformly to a single rank $n+1$ projection $\Pi_\infty$. Every subsequential hemisphere therefore lies in the same plane $V=\operatorname{im}\Pi_\infty$. The boundary condition gives $q/R\in V$, so the centre, radius, plane and side of the limiting hemisphere are all fixed.

After flattening the rescaled support sphere, boundary preserving graph parametrisations over this fixed hemisphere identify the submanifolds for all times sufficiently close to $T$ with the same domain. Every sequence has a subsequence converging smoothly to the inclusion in these parametrisations. Thus the entire family
	\begin{align*}
	\frac{F(\cdot,t)-q}{\sqrt{2n(T-t)}}
	\end{align*}
converges smoothly, up to the free boundary, to the unit hemisphere in $V$.
\end{proof}

\appendix
\section{Compactness and continuation for spherical free boundary flow}\label{app:fb-compactness}

The compactness theory for immersed submanifolds originates with Langer's theorem for surfaces \cite{Langer1985} and was extended to arbitrary dimension and codimension by Breuning \cite{Breuning2015}. We adapt this graphical approach to spherical free boundary flow, retaining the boundary regularity, reflection and continuation arguments needed in Section \ref{sec:fb-convergence}.

\begin{lemma}[Smooth boundary compactness and continuation]
\label{lem_fb-smooth-compactness}
Let $F_j:M_j^n\times[a,b]\to\overline B_{R_j}(c_j)$ be smooth compact connected spherical free boundary flows, where $n,m\ge1$, $a<b$ and $R_j\geq R_*>0$, satisfying
	\begin{align}\label{eq_fb-compactness-input}
	|A_j|&\leq K_0,\qquad |\nabla A_j|\leq K_1,\qquad |\nabla^2A_j|\leq K_2.
	\end{align}
Then:
\begin{enumerate}[(i)]
\item For every $k\ge0$ and $0<\delta<b-a$,
	\begin{align*}
	\sup_{M_j\times[a+\delta,b]}|\nabla^kA_j|&\leq C_k,
	\end{align*}
including the boundary, with constants depending only on $n,m,R_*,K_0,K_1,K_2,k,\delta$.
\item Given $t_*\in(a,b]$ and points $p_j$ with $F_j(p_j,t_*)$ bounded, a subsequence converges smoothly on compact subsets of space and on every $[a+\delta,b]$ to a complete pointed flow. Convergence is through time independent boundary preserving embeddings whose images exhaust intrinsic balls about $p_j$ at $t_*$. Completeness includes the boundary.

The limit has empty boundary when the intrinsic distances of $p_j$ from $\partial M_j$ diverge, and nonempty boundary when they remain bounded. In the latter case, if $R_j\to\infty$, the limiting support is a hyperplane and reflection gives a smooth complete flow without boundary.
\item A compact spherical free boundary flow on $[0,T)$, with $T<\infty$, extends smoothly past $T$ if the above three curvature quantities remain bounded.
\end{enumerate}
\end{lemma}

\begin{proof}
We use Breuning's local graph construction \cite[Section 2]{Breuning2015} and take diagonal subsequences on intrinsic balls. The required boundary charts follow by making the ambient metric a product near the support and doubling the immersion. Orthogonality gives a $C^{1,1}$ double with bounded second fundamental form, so projection onto tangent planes supplies uniform graphical charts over balls and half balls.

In Fermi coordinates, with $s$ the inward distance from the support, the graphical equation and boundary condition are
	\begin{align*}
	\partial_tu^\alpha&=g^{ij}(y,s,u,Du)\partial_{ij}u^\alpha+b^\alpha(y,s,u,Du),\\
	\partial_su^\alpha&=0,\qquad\text{ on }s=0.
	\end{align*}
The principal operator is uniformly parabolic and the boundary operator is componentwise Neumann. Thus the complementing condition holds in arbitrary codimension. The assumed curvature bounds and boundary Schauder estimates \cite{Solonnikov1965} give (i); pointed compactness in these charts gives (ii). When the limiting support is a hyperplane, even reflection of the Neumann graphs and parabolic regularity give the asserted smooth double.

For (iii), the higher derivative bounds give a smooth limiting immersion $F_T$, satisfying the support and angle conditions and all compatibility conditions. In normal graph coordinates that preserve the boundary position, the linearised angle operator is
	\begin{align*}
	B_0V&=\nabla_\mu^\perp V-R^{-1}V.
	\end{align*}
Its principal part is again Neumann. The standard local existence argument for quasilinear equations based on the linear solvability theorem \cite[Theorem 4.9]{Solonnikov1965} therefore restarts the flow from $F_T$; compatibility gives smoothness across $T$. A boundary preserving reparametrisation restores the normal velocity equation.
\end{proof}


\begin{thebibliography}{99}

\bibitem{AndrewsBaker}
B. Andrews and C. Baker.
\newblock \href{https://doi.org/10.4310/jdg/1292940688}{Mean curvature flow of pinched submanifolds to spheres}.
\newblock \emph{J. Differential Geom.} \textbf{85} (2010), no.3, 357--395.

\bibitem{Breuning2015}
P. Breuning.
\newblock \href{https://doi.org/10.1007/s12220-014-9472-7}{Immersions with bounded second fundamental form}.
\newblock \emph{J. Geom. Anal.} \textbf{25} (2015), no.2, 1344--1386.

\bibitem{Buckland}
J. A. Buckland.
\newblock \href{https://doi.org/10.1515/crll.2005.2005.586.71}{Mean curvature flow with free boundary on smooth hypersurfaces}.
\newblock \emph{J. Reine Angew. Math.} \textbf{586} (2005), 71--90.

\bibitem{Chow2026}
T.-K. A. Chow.
\newblock Tensorial pinching by deforming curvature cones in higher codimension mean curvature flow.
\newblock Preprint (2026), \href{https://arxiv.org/abs/2610.25610}{\nolinkurl{arXiv:2610.25610}}.

\bibitem{EdelenConvexity2016}
N. Edelen.
\newblock \href{https://doi.org/10.1016/j.aim.2016.02.026}{Convexity estimates for mean curvature flow with free boundary}.
\newblock \emph{Adv. Math.} \textbf{294} (2016), 1--36.

\bibitem{Edelen2016}
N. Edelen.
\newblock \href{https://doi.org/10.1515/crelle-2017-0053}{The free-boundary Brakke flow}.
\newblock \emph{J. Reine Angew. Math.} \textbf{758} (2020), 95--137.

\bibitem{Edelenetal}
N. Edelen, R. Haslhofer, M. N. Ivaki and J. J. Zhu.
\newblock \href{https://doi.org/10.1002/cpa.22009}{Mean convex mean curvature flow with free boundary}.
\newblock \emph{Comm. Pure Appl. Math.} \textbf{75} (2022), no.4, 767--817.

\bibitem{EvansLambertWood2022}
C. G. Evans, B. Lambert and A. Wood.
\newblock \href{https://doi.org/10.1007/s00526-022-02229-0}{Lagrangian mean curvature flow with boundary}.
\newblock \emph{Calc. Var. Partial Differential Equations} \textbf{61} (2022), article no.106.

\bibitem{HamiltonCompactness1994}
R. S. Hamilton.
\newblock \href{https://doi.org/10.4310/CAG.1994.v2.n1.a10}{Convex hypersurfaces with pinched second fundamental form}.
\newblock \emph{Comm. Anal. Geom.} \textbf{2} (1994), no.1, 167--172.

\bibitem{Hamilton}
R. S. Hamilton.
\newblock \href{https://doi.org/10.4310/jdg/1214456010}{Harnack estimate for the mean curvature flow}.
\newblock \emph{J. Differential Geom.} \textbf{41} (1995), no.1, 215--226.

\bibitem{Huisken1984}
G. Huisken.
\newblock \href{https://doi.org/10.4310/jdg/1214438998}{Flow by mean curvature of convex surfaces into spheres}.
\newblock \emph{J. Differential Geom.} \textbf{20} (1984), 237--266.

\bibitem{HuiskenSinestrari2015}
G. Huisken and C. Sinestrari.
\newblock \href{https://doi.org/10.4310/jdg/1442364652}{Convex ancient solutions of the mean curvature flow}.
\newblock \emph{J. Differential Geom.} \textbf{101} (2015), no.2, 267--287.

\bibitem{Langer1985}
J. Langer.
\newblock \href{https://doi.org/10.1007/BF01456183}{A compactness theorem for surfaces with $L_p$-bounded second fundamental form}.
\newblock \emph{Math. Ann.} \textbf{270} (1985), no.2, 223--234.

\bibitem{Lyn}
S. Lynch.
\newblock The planarity estimate in dimensions five and six.
\newblock Preprint (2026), \href{https://arxiv.org/abs/2609.39330v1}{\nolinkurl{arXiv:2609.39330v1}}.

\bibitem{LNAncient2021}
S. Lynch and H. T. Nguyen.
\newblock \href{https://doi.org/10.1007/s00526-020-01888-1}{Pinched ancient solutions to the high codimension mean curvature flow}.
\newblock \emph{Calc. Var. Partial Differential Equations} \textbf{60} (2021), Article29.

\bibitem{LN}
S. Lynch and H. T. Nguyen.
\newblock \href{https://doi.org/10.1007/s00208-022-02545-y}{Convexity estimates for high codimension mean curvature flow}.
\newblock \emph{Math. Ann.} \textbf{388} (2024), 575--613.

\bibitem{LNSurgery2025}
S. Lynch and H. T. Nguyen.
\newblock High codimension mean curvature flow with surgery.
\newblock Preprint (2025), \href{https://arxiv.org/abs/2004.07163v3}{\nolinkurl{arXiv:2004.07163v3}}.

\bibitem{MichaelSimon1973}
J. H. Michael and L. M. Simon.
\newblock \href{https://doi.org/10.1002/cpa.3160260305}{Sobolev and mean-value inequalities on generalized submanifolds of $\mathbb R^n$}.
\newblock \emph{Comm. Pure Appl. Math.} \textbf{26} (1973), 361--379.

\bibitem{Naff2022}
K. Naff.
\newblock \href{https://doi.org/10.1215/00127094-2021-0026}{A planarity estimate for pinched solutions of mean curvature flow}.
\newblock \emph{Duke Math. J.} \textbf{171} (2022), no.2, 443--482.

\bibitem{FBHCCSF}
H. T. Nguyen and A. A. Vogiatzi.
\newblock High Codimension Curve Shortening Flow with Free Boundary.
\newblock Preprint (2026), \href{https://arxiv.org/abs/2602.20865}{\nolinkurl{arXiv:2602.20865}}

\bibitem{Solonnikov1965}
V. A. Solonnikov.
\newblock On boundary value problems for linear parabolic systems of differential equations of general form.
\newblock \emph{Proc. Steklov Inst. Math.} \textbf{83} (1965), 1--184.
\newblock \url{https://www.mathnet.ru/eng/tm2755}.

\bibitem{Stahl1996}
A. Stahl.
\newblock \href{https://doi.org/10.1007/BF01190825}{Regularity estimates for solutions to the mean curvature flow with a Neumann boundary condition}.
\newblock \emph{Calc. Var. Partial Differential Equations} \textbf{4} (1996), no.4, 385--407.

\bibitem{StahlConvergence1996}
A. Stahl.
\newblock \href{https://doi.org/10.1007/BF01246150}{Convergence of solutions to the mean curvature flow with a Neumann boundary condition}.
\newblock \emph{Calc. Var. Partial Differential Equations} \textbf{4} (1996), no.5, 421--441.

\end{thebibliography}
\end{document}